\documentclass[12pt]{amsart}

\usepackage{graphicx}
\usepackage{amscd}
\usepackage{amsmath}
\usepackage{amsfonts}
\usepackage{amssymb}
\usepackage{a4wide}
\usepackage{setspace}
\usepackage{enumerate}         
\usepackage{color}
\usepackage{url}
\usepackage{amsthm}
\usepackage{hyperref}
\usepackage{bm}
\usepackage{xy}
\usepackage{stmaryrd}

\theoremstyle{plain}
\newtheorem{theorem}{Theorem}[section]

\newtheorem{lemma}[theorem]{Lemma}

\newtheorem{proposition}[theorem]{Proposition}

\newtheorem{definition}[theorem]{Definition}

\newtheorem{assumption}[theorem]{Assumption}
\theoremstyle{remark}
\newtheorem{remark}[theorem]{Remark}

\numberwithin{equation}{section}

\newcommand{\ind}{1\!\kern-1pt \mathrm{I}}
\newcommand{\rsto}{]\!\kern-1.8pt ]}
\newcommand{\lsto}{[\!\kern-1.7pt [}

\numberwithin{equation}{section}

\begin{document}
\title[Quadratic BSDEs driven by G-Brownian motion] {Quadratic backward stochastic differential equations driven by G-Brownian motion: discrete solutions and approximation}
\thanks{Corresponding author: yiqing.lin@polytechnique.edu (Y. Lin)}
\author{Ying Hu}
\address{Ying Hu\newline\indent Institut de recherche math\'ematiques de Rennes
\newline\indent Universit\'e de Rennes 1\newline\indent
35042 Rennes Cedex - France}

\author{Yiqing Lin}
\address{Yiqing Lin\newline\indent Centre de math\'ematiques appliquées \newline\indent \'Ecole Polytechnique\newline\indent 91128 Palaiseau Cedex - France}

\author{Abdoulaye Soumana Hima}
\address{Abdoulaye Soumana Hima \newline\indent Institut de recherche math\'ematiques de Rennes
\newline\indent Universit\'e de Rennes 1\newline\indent
35042 Rennes Cedex - France \newline\indent and\newline\indent 
D\'epartement de math\'ematiques  \newline\indent 
Universit\'e de Maradi\newline\indent 
BP 465 Maradi - Niger}

\begin{abstract}
In this paper, we consider backward stochastic differential equations driven by $G$-Brownian motion (GBSDEs) under quadratic assumptions on coefficients. We prove the existence and uniqueness of solution for such equations. On the one hand, a priori estimates are obtained by applying the Girsanov type theorem in the $G$-framework, from which we deduce the uniqueness. On the other hand, to prove the existence of solutions, we first construct solutions for discrete GBSDEs by solving corresponding fully nonlinear PDEs, and then approximate solutions for general quadratic GBSDEs in Banach spaces.
\end{abstract}

\date{\today}
\maketitle


\section{Introduction}
The first existence and uniqueness result for nonlinear backward stochastic differential equations (BSDEs for short) of the following form is provided by Pardoux and Peng in \cite{PP}:
\begin{equation*}
Y_{t}=\xi +\int_{t}^{T}f\left( s,Y_{s},Z_{s}\right)
ds-\int_{t}^{T}Z_{s}dW_s,\ 0\leq t\leq T,
\end{equation*}%
where the generator $f$ is uniformly Lipschitz  and the terminal value $\xi $ is square integrable.
 Since then, BSDEs have been
studied with great interest and moreover, these equations are found to have 
strong 
connections with different mathematical fields, such as mathematical
finance, stochastic control and partial differential equations. In
particular, many efforts have been made to relax the assumption on the
generator. For instance, Kobylanski \cite{Kob} was the first to investigate the BSDE with a generator having quadratic growth in $Z$ and a bounded terminal value. She used an exponential transformation in order to come back to the framework of linear growth generator. This seminal work of quadratic BSDE has been extended by many authors. Since a complete review of these literatures is too extensive, we only concentrate on those of immediate interest. For the existence, Briand and Hu \cite{BH06} observed that the existence of exponential moments of the terminal condition is sufficient to construct a solution of quadratic BSDEs; uniqueness is proved in \cite{BH08} under the additional assumption that the generator is convex (or concave). On the other hand, Hu et al. discovered in \cite{HIM} that  for certain type of $f$ locally Lipschitz in $Z$, if the solution $Y$ is bounded, the solution $Z$ is  bounded in BMO norm, and thus the uniqueness could be proved by applying linearization of the generator. Afterwards, Tevzadze  improved the methodology of \cite{HIM}  and he gave a direct proof in \cite{Tev} for the solvability of quadratic BSDEs by standard fixed point arguments, whereas the terminal value was assumed to be sufficient small. To get rid of this technical  assumption on $\xi$, Briand and Elie exhibited a priori estimates in the light of \cite{AIDR, B08} and approximated bounded terminal values with Malliavin differentiable ones in \cite{BE}.  In more general situations, \
Morlais reconsidered the problem of \cite{Kob} with continuous martingale driver in  \cite{Mor}, while recently Barrieu and El Karoui obtained similar results in  \cite{BEL} but by a completely different forward method.\\
%

Motivated by mathematical finance problems with Knightian uncertainty, Peng established systematically in  \cite{Peng2007, Peng2008, Pen10, Peng2011} a framework of time-consistent sublinear expectation, called $G$-expectation. In particular, this sublinear expectation is associated with a new type of Brownian motion $(B_t)_{t\geq 0}$, i.e., $G$-Brownian motion, which 
has independent, stationary and $G$-normally distributed increments. This process and its quadratic variation $\langle B\rangle$ play center roles in the related nonlinear stochastic analysis.
Indeed, the stochastic integrals with respect to $G$-Brownian motion and its quadratic variation have been first introduced by Peng in his pioneer work \cite{Peng2007}, which are initially defined on the simple process space and later extended as linear operator on  Banach completions. Thereafter, the $G$-stochastic calculus is further developed, for example, in \cite{Pen10, G, LP11, Lin}. Another important feature of the $G$-expectation is found by Denis et al. in \cite{DHP11}, namely,  the $G$-expectation can be represented by the upper expectation over a collection of mutually singular martingale measures $\mathcal{P}_G$. Moreover, the notion of quasi-sure with respect to the associated Choquet capacity is introduced by Denis et al. to the $G$-framework. \\

As their classical counterparts, stochastic differential equations driven by $G$-Brownian motion (GSDEs) are well defined in the quasi-sure sense and their solvability can be established by the contracting mapping theorem under Lipschitz assumptions (cf. \cite{Pen10} and \cite{G}). However, the challenging problem of wellposedness for backward GSDEs (GBSDEs)  remained open until a complete theorem has been proved by Hu et al. \cite{HJPS14a}.\\

Similarly to the classical case, the $G$-martingale representation theorem is heuristic to the formulation of GBSDEs, which reads as follows
\begin{equation}\label{mar1}
M_t=M_0+ \overline{M}_t+K_t,
\end{equation}
where 
\begin{equation}\label{mar2}
\overline{M}_t=\int^t_0 Z_s dB_s\quad\mbox{and}\quad K_t=\int^t_0 \eta_s d\langle B\rangle_s-\int^t_0 2G(\eta_s)ds.
\end{equation}
In contrast to the classical martingale representation, the $G$-martingale $M$ is decomposed into two parts: the $G$-It\^o type integral part $\overline{M}=\int Z dB$, which is called symmetric $G$-martingale, in the sense that $-\overline{M}$ is still a $G$-martingale; the decreasing $G$-martingale part $K$, which vanishes in the classical theory, however, plays a significant role in this new context.
Whether the process $K$ admits a unique representation in the form (\ref{mar2}) is a sophisticated question. The first positive answer is given by Peng in \cite{Pen07a} for the $G$-martingale associated with a terminal value $M_T\in Lip(\Omega_T)$, which reads as smooth and finitely dimensional path function.  It is also worth mentioning that a series of successive works by Soner et al. \cite{STZ11} and Song \cite{song2011some} affirm the existence and uniqueness of the first level decomposition (\ref{mar1}) for $M_T\in L^p_G(\Omega)$, $p>1$, which is the Banach completion of $Lip(\Omega_T)$. Finally, with the help of the norm creatively introduced in Song \cite{song2012uniqueness}, a complete theorem for $G$-martingale  representation has been obtained by Peng et al. \cite{PSZ} on a  complete metric subspace of $L^p_G(\Omega)$, $p>1$.\\ 

We take into consideration of the $G$-martingale representation theorem and naturally, we can formulate GBSDE as follows, where the decreasing $G$-martingale $K$ appears in the dynamics:
\begin{equation}\label{gb1}
Y_t=\xi +\int^T_t g(s, Y_s, Z_s) ds+\int^T_t h(s, Y_s, Z_s) d\langle B\rangle_s -\int^T_t Z_s dB_s-(K_T-K_t). 
\end{equation}
Under Lipschitz assumptions on the generators, Hu et al. investigated in \cite{HJPS14a} the existence and uniqueness of the triple $(Y, Z, K)$ in proper Banach spaces satisfying the above equation. They started with BSDEs with bounded and smooth generators and Markovian terminal values and constructed solutions by classical solutions  of fully nonlinear PDEs (cf. Krylov's results in \cite{Kry87}). Then, the partition of unity theorem was employed in \cite{HJPS14a} to proceed a type of Galerkin approximation to solutions of GBSDEs with general  parameters. Besides, the uniqueness was deduced in \cite{HJPS14a} based on a priori estimates. In particular, the uniqueness of  $K$ is impressive in the light of $G$-martingale estimates found in \cite{song2011some}. The results in \cite{HJPS14a} breaks new ground in the $G$-expectation theory. In the accompanying paper \cite{HJPS14b}, Hu et al. discussed fundamental properties of the above GBSDE: the comparison theorem, the fully nonlinear Feynman-Kac formula and the related Girsanov transformation. Moreover, the correspondence between GBSDEs and Sobolev type solutions of nonlinear path-dependent PDEs is examined in \cite{PS}.\\

We now compare the result of \cite{HJPS14a} with the profound works \cite{STZ12, STZ13} by Soner et al.,  in which the so-called second order backward stochastic differential equations (2BSDEs) are deeply studied. 
 This type of equation is highly related to the GBSDE and it is defined on the Wiener space as follows:
$$
Y_t=\xi +\int^T_t F(s, Y_s, Z_s, \widehat{a}_s) ds -\int^T_t Z_s dB_s+(K_T-K_t),\quad \mathbb{P}-\mbox{a.s.},\ \mbox{for all}\ \mathbb{P}\in \mathcal{P}_H, 
$$
where $B$ is the canonical process, the process $\widehat{a}$ is the density of $\langle B\rangle$ and $\mathcal{P}_H$ is a collection of martingale measures similar to $\mathcal{P}_G$ (could be even larger).  This equation is a reinforced BSDE in the sense that it holds true $\mathbb{P}$-a.s. for all $\mathbb{P}\in \mathcal{P}_H$ and moreover, the family of $K:=\{K^\mathbb{P}\}_{\mathbb{P}\in \mathcal{P}_H}$ should satisfy a minimum condition (then $-K$ verifies the $G$-martingale constraint in the GBSDE context, see \cite{STZ11}):
$$
K^\mathbb{P}_t=\mathop{\mbox{essinf}^\mathbb{P}}\limits_{\mathbb{P}'\in \mathcal{P}_H(t, \mathbb{P})} {E}^{\mathbb{P}'}_t[K^{\mathbb{P}'}_T],\quad \mathbb{P}-\mbox{a.s.},\ \mbox{for all}\ \mathbb{P}\in \mathcal{P}_H, \quad 0\leq t\leq T.
$$
Under Lipschitz assumptions, the uniqueness of the 2BSDE is proved in \cite{STZ12} by observing that the solution to the 2BSDE can be represented as the (essential) supremum of a class of martingale-driven BSDEs solutions. For the existence, the proof involves a delicate pathwise construction: the process $Y$ is defined pathwisely by solutions of BSDEs on shift spaces. This process verifies a critical principle of optimality and thus, the structure of 2BSDE could be derived from the $g$-supermartingale decomposition (cf. \cite{Pen99}), where the family of processes $K$ can be a posteriori aggregated once the stochastic integral part is aggregated by Nutz \cite{N12}. To get rid of the measurability problem during the construction of solutions, Soner et al. assume the technical condition that both $\xi$ and $F$ is uniformly continuous in $\omega$, whereas this assumption is
removed in the recent work of Possama\"i et al. \cite{PTZ}. In the framework of 2BSDEs, the results in \cite{STZ12, STZ13} are generalized by Possama\"i and Zhou \cite{PZ13} and by Lin \cite{Lin16} to the quadratic case and furthermore, Matoussi et al. \cite{MPZ} applied quadratic 2BSDEs to solve the utility maximization problems from \cite{HIM} in the context with non-dominated models.  One could see that the GBSDE (\ref{gb1}) actually corresponds to  the 2BSDE defined with
$$
F(t, y, z, a)=g(t, y, z)+h(t, y, z)a,
$$
however, the GBSDE requires more structure conditions on the coefficient and the terminal value so that the solution can be found with more regularity adapted to the requirement of process space in the $G$-framework.\\

The main objective of this paper is to provide the existence and uniqueness result for scalar-valued quadratic GBSDEs adapted to the setting of \cite{HJPS14a}. Without loss of generality, we consider only the following type of GBSDE:
\begin{equation*}
Y_t=\xi +\int^T_t h(s, B_{\cdot\wedge s}(\cdot), Y_s, Z_s) d\langle B\rangle_s -\int^T_t Z_s dB_s-(K_T-K_t),
\end{equation*}
where $\xi$ is an element in the $L^\infty_G$ completion of $Lip(\Omega_T)$ and $h$ is Lipschitz in $y$ and locally Lipschitz in $z$, similarly to \cite{HIM} and \cite{Mor} in the classical framework. Moreover, we require that $h$ is uniformly continuous with respect to $\omega$, which is a technical condition for the successive approximation. This assumption is stronger than the corresponding structure condition in \cite{HJPS14a} and how to weaken such technical assumption is postponed to future research.\\

   This paper is organized as follows. Section 2 is dedicated to preliminaries in the $G$-framework and the formulation of quadratic GBSDEs. In Section 3, we introduce the space of $G$-BMO martingale generators and deduce a priori estimates for quadratic GBSDEs through the $G$-Girsanov transformation. Meanwhile, we obtain the uniqueness straightforwardly by the a priori estimates.  In Section 4, we consider GBSDEs with discrete generators and terminal values of the following functional type:
   \begin{eqnarray*}
\xi  &=&\varphi \left( B_{t_{1}},B_{t_{2}}-B_{t_{1}},\ldots
,B_{t_{N}}-B_{t_{N-1}}\right);  \\
h\left( t,B_{\cdot\wedge t}(\cdot),Y_{t},Z_{t}\right)  &=&f\left( t,B_{t_{1}\wedge
t},B_{t_{2}\wedge t}-B_{t_{1}\wedge t},\ldots ,B_{t_{N}\wedge
t}-B_{t_{N-1}\wedge t},Y_{t},Z_{t}\right),
\end{eqnarray*}%
where $0=t_{0}\leq t_{1}\leq \ldots \leq t_{N-1}\leq t_{N}=T$ is a given partition of $\left[
0,T\right]$. We proceed with the argument of Hu and Ma in \cite{HM04} to construct solutions of such GBSDEs, where Krylov's estimates for fully nonlinear PDEs are applied  as in \cite{HJPS14a}. The last section present the technics of discretization and regularization, moreover,  
the existence result for general quadratic GBSDEs shall then
be proved  by successive approximation.

\section{Preliminary}
\subsection{The $G$-framework}
In this section, we review notations and basic results in the framework of $G$-expectation, which concern the formulation of $G$-Brownian motion and related $G$-stochastic calculus. In this paper, we only consider the one-dimensional case. The readers interested in more details on this topic are referred to \cite{Pen10, G, LP11, Lin}. \\

Let  $\Omega$ be a complete separable metric space, and let 
 $\mathcal{H}$ be a linear
space of real-valued functions defined on $\Omega$ satisfying: if
$X_{i}\in \mathcal{H}$, $i=1,\ldots,n$,
then%
\[
\varphi(X_{1}, X_2, \ldots,X_{n})\in \mathcal{H}, \quad {\forall }\varphi \in
\mathcal C_{l,Lip}(\mathbb{R}^{n}),
\]
where $\mathcal C_{l,Lip}(\mathbb{R}^{n})$ is the space of all
continuous real-valued functions defined on $\mathbb{R}^n$ such that for some $C>0$ and $k\in \mathbb{N}$ depending on $\varphi$,
$$|\varphi(x)-\varphi(y)|\leq C(1+|x|^k+|y|^k)|x-y|,\quad \forall
x,y\in\mathbb{R}^n.$$


\begin{definition}[Sublinear expectation space]
A sublinear expectation $\mathbb{E}[\cdot]$ is a functional 
$\mathbb{E}:\mathcal{H}\rightarrow \mathbb{R}$ satisfying the following
properties: for all $X$, $Y\in \mathcal{H}$, we have\newline
(a) Monotonicity: if $X\geq Y$, then $\mathbb{E}[X]\geq \mathbb{E}[Y]$;%
\newline
(b) Constant preservation: $\mathbb{E}[c]=c$, $c\in \mathbb R$;\newline
(c) Sub-additivity: $\mathbb{E}[X+Y]\leq \mathbb{E}[X]+\mathbb{E}[Y]$;%
\newline
(d) Positive homogeneity: $\mathbb{E}[\lambda X]=\lambda \mathbb{E}[X]$, for
all $\lambda \geq 0$.\\
We call the triple $(\Omega ,\mathcal{H},\mathbb{E})$ sublinear
expectation space.
\end{definition}
\begin{definition}[Independence]
Fix the sublinear expectation space $(\Omega,\mathcal{H},\mathbb{E})$.  A random
variable $Y\in \mathcal{H}$
 is said to be independent of $X_1,  X_2, \ldots, X_{n}\in \mathcal{H}$, if for all $\varphi \in
\mathcal C_{l,Lip}(\mathbb{R}^{n+1})$,
\[
\mathbb{{E}}\left[\varphi(X_1, X_2, \ldots, X_{n},Y)\right]=\mathbb{{E}}\left[\mathbb{{E}}%
\left[\varphi(x_1, x_2, \ldots, x_{n},Y)\right]\big|_{(x_1, x_2, \ldots, x_{n})=(X_1, X_2, \ldots, X_{n})}\right].
\]
\end{definition}

Now we introduce the definition of $G$-normal distribution. 
\begin{definition}[$G$-normal distribution]
We say the random variable $X\in \mathcal{H}$ 
is $G$-normally distributed, noted by $X\sim \mathcal{N}(0, [\underline{\sigma}, \overline{\sigma}])$, $0\leq \underline{\sigma}\leq \overline{\sigma}$, if for
any function $\varphi \in \mathcal C_{l.Lip}(\mathbb{R})$, the function $u$  
defined by $u(t,x):=\mathbb{E}[\varphi (x+\sqrt{t}X)]$, $\left( t,x\right) \in %
\left[ 0,\infty \right) \times \mathbb{R}$, is a viscosity solution of
$G$-heat equation:
\[
\partial _{t}u-G\left(D_{x}^{2}u\right)=0;~~u(0,x)=\varphi (x), 
\]%
where \[
G(a):=\frac{1}{2}\left(\overline{\sigma }%
^{2}a^{+}-\underline{\sigma }^{2}a^{-}\right).
\]
%
\end{definition}
Throughout this paper, we consider only the non-degenerate case, i.e., $\underline{\sigma}>0$.
We now fix $\Omega:=\mathcal{C}[0, \infty)$, which is equipped with the raw filtration $\mathcal{F}$ generated by the canonical process $(B_t)_{t\geq0}$,
i.e., $B_t(\omega)=\omega_t$, for $(t, \omega) \in[0,\infty)\times\Omega$. Let us consider the function spaces defined by
\begin{align*}Lip(\Omega _T):=\{\varphi(B_{t_1},B_{t_2}-B_{t_1},&\ldots,B_{t_n}-B_{t_{n-1}}): \\
&\ 0\leq t_{1}\leq t_2 \ldots ,t_{n}\leq T,\ \varphi\in \mathcal C_{l,Lip}(\mathbb{R}^n)\},\quad \mbox{for}\ T>0,
\end{align*}
and
$$Lip(\Omega):=\bigcup_{n=1}^\infty Lip(\Omega_n).$$

\begin{definition}[$G$-Brownian motion and $G$-expectation]
 On the sublinear expectation space $\left( \Omega , Lip(\Omega),\mathbb{E}\right)$, the canonical process $\left(B_{t}\right) _{t\geq 0}$ 
 is called $G$-Brownian motion if the following properties are verified:\newline
(a) $B_{0}=0;$\newline
(b) for each $t$, $s\geq 0$, the increment $B_{t+s}-B_{t}\sim \mathcal{N}(0, [\sqrt{s}\underline{\sigma}, \sqrt{s}\overline{\sigma}])$
 and is independent from $\left( B_{t_{1}},B_{t_{2}},\ldots
,B_{t_{n}}\right) $, for $0\leq t_{1}\leq t_{2}\leq
\ldots \leq t_{n}\leq t$.\newline
Moreover, the sublinear expectation $\mathbb{E}[\cdot]$ is called $G$-expectation.
\end{definition}
%

\begin{definition}[Conditional $G$-expectation]
For the random variable $\xi\in Lip(\Omega_T)$ of the following form:
$$\varphi(B_{t_1},B_{t_2}-B_{t_1},\ldots,B_{t_n}-B_{t_{n-1}}),\quad \varphi\in \mathcal C_{l,Lip}(\mathbb{R}^n), $$
the conditional $G$-expectation $\mathbb{E}_{t_i}[\cdot]$, $i=1, \ldots, n$, is defined as follows
\[
\mathbb{E}_{t_{i}}[\varphi (B_{t_{1}}, B_{t_{2}}-B_{t_{1}},\ldots
,B_{t_{n}}-B_{t_{n-1}})]=\widetilde{\varphi}%
(B_{t_{1}}, B_{t_{2}}-B_{t_{1}},\ldots ,B_{t_{i}}-B_{t_{i-1}}),
\]
where %
\[
\widetilde{\varphi }\left( x_{1},\ldots ,x_{i}\right) =\mathbb{E}\left[
\varphi \left( x_{1},\ldots ,x_{i},B_{t_{i+1}}-B_{t_{i}},\ldots
,B_{t_{n}}-B_{t_{n-1}}\right) \right] .
\]
If $t\in (t_i, t_{i+1})$, then the conditional $G$-expectation $\mathbb{E}_t[\xi]$ could be defined by reformulating $\xi$ as
\begin{align*}
\xi=\widehat{\varphi}(B_{t_1},B_{t_2}-B_{t_1},\ldots,B_t-B_{t_i}, B_{i+1}-B_t, \ldots, B_{t_n}-B_{t_{n-1}}), \quad \widehat\varphi\in \mathcal C_{l,Lip}(\mathbb{R}^{n+1})
\end{align*}
\end{definition}
For $\xi \in Lip(\Omega_T)$ and $p\geq 1$, we consider the norm $\left\Vert \xi
\right\Vert _{L^p_G}=(\mathbb{E}[|\xi |^{p}])^{1/p}$. Denote by $L^p_G(\Omega_T)$ the Banach completion of $Lip(\Omega_T)$ under $\Vert\cdot\Vert_{L^p_G}$. It is easy to check that the conditional $G$-expectation $\mathbb{E}_t[\cdot]: Lip(\Omega_T)\rightarrow Lip(\Omega_t)$ is a continuous mapping and thus can be extended to $\mathbb{E}_t: L^p_G(\Omega_T)\rightarrow L^p_G(\Omega_t)$.
\begin{definition}[$G$-martingale]
A process $\left( M_{t}\right) _{t\in \left[ 0,T\right] }$ with $M_t\in 
L_{G}^{1}(\Omega _{t})$, $0\leq t\leq T$, is a $G$-martingale if $\mathbb{E}_{s}[M_{t}]=M_{s}$,
for all $0\leq s\leq t\leq T$. The process $\left( M_{t}\right) _{t\in \left[ 0,T\right]}$ is called symmetric $G$-martingale if $-M$ is also a $G$-martingale. 
\end{definition}

According to Denis et al. \cite{DHP11}, we have the following representation theorem of $G$-expectation on $L^1_G(\Omega_T)$. In the sequel, we denote by $\mathbb P_0$ the Wiener measure, under which the canonical process $(B_t)_{t\geq 0}$ is a $\mathbb P_0$-Brownian motion.
\begin{theorem} [Representation of $G$-expectation]The $G$-expectation can be represented by the upper expectation over a collection of martingale measures, i.e., for $\xi\in L^1_G(\Omega_T)$, we have
\[
\mathbb{E}[\xi ]=\sup_{\mathbb P\in \mathcal{P}_G}E^{\mathbb P}[\xi],
\]
where 
\begin{align*}
\mathcal{P}_{G}=\bigg\{ \mathbb P_{h}: \mathbb P_{h}=\mathbb P_{0}\circ
X^{-1},~X_{t}=\int_{0}^{t}h_{s}dB_{s},&~h\in \mathbb{H}^2_{\mathbb P_0}(0, T),\\
 &h_t\in [\underline{\sigma },\overline{\sigma }]),\ \mathbb P_0-a.s.,\ 0\leq t\leq T\bigg\}.
\end{align*}
\end{theorem}

By the theorem above, the $G$-expectation can be extended to a larger domain, i.e., for all $\mathcal{F}_T$ measurable function $X$, $\mathbb{E}[X]:=\sup_{\mathbb P\in \mathcal{P}_G}E^{\mathbb P}[X]$. It is also proved in \cite{DHP11} that $\mathcal{P}_G$ is relatively weakly compact and thus its completion $\overline{P}_G$ is weakly compact. Therefore, we can naturally define the Choquet capacity $\overline C(\cdot)$ by $\overline{C}(A):=\sup_{P\in \overline{\mathcal P}_G} P(A)$, $A\in \mathcal{B}(\Omega_T)$ and introduce the notion of quasi-sure.

\begin{definition}[Quasi-sure]
A set $A\in \mathcal{B}(\Omega_{T})$ is a polar if $\overline{C}(A)=0$. A property holds
``quasi-surely" (q.s.) if it is true outside a polar set.
\end{definition}
The following proposition helps to understand the correspondence between GBSDEs and 2BSDEs.
\begin{proposition}\label{emu} Suppose $X$ and $Y\in L^1_G$, then the following statements are equivalent:\newline
(a) for each $\mathbb{P}\in \mathcal{P}_G$, $X=Y$, $\mathbb{P}$-a.s.;\\
(b) $\mathbb{E}[|X-Y|]=0$;\\
(c) $\overline{C}(\{X\neq Y\})=0$.
\end{proposition}
For the terminal value of quadratic GBSDE, we define the space $L^\infty_G(\Omega_T)$ as the completion of $Lip(\Omega_T)$ under the norm  
$$
\Vert\xi\Vert_{L^\infty_G}:=\inf\left\{M\geq 0: |\xi|\leq M\ q.s.\right\}.
$$

A important feature of the $G$-expectation theory is that the quadratic variation of the $G$-Brownian motion  $(\langle B\rangle_t)_{t\geq 0}$ is no longer a deterministic process, which
 is given by
\[
\left\langle B\right\rangle _{t}=\lim_{\mu \left( \pi _{t}^{N}\right)
\rightarrow
0}\sum_{j=0}^{N-1}(B_{t_{j+1}^{N}}-B_{t_{j}^{N}})^{2}=B_{t}^{2}-2%
\int_{0}^{t}B_{s}dB_{s},
\]%
where $\pi _{t}^{N}$, $N=1,2,\mathbb{\ldots }$, are refining partitions of $[0, t]$.
By Peng \cite{Pen10}, for all $t$, $s\geq 0$,
$\langle B\rangle_{t+s}-\langle B\rangle_{t}\in \left[s\underline{\sigma}^{2}, s%
\overline{\sigma}^{2}\right]$, q.s.\\

In what follows, we discuss the stochastic integrals with respect to the $G$-Brownian motion and its quadratic variation. 
\begin{definition}
Let $\mathcal H_{G}^{0}\left( 0,T\right) $ be the set of simple processes of the following form:
\begin{equation}\label{simp}
\eta _{t}\left( \omega \right) =\sum_{j=0}^{N-1}\xi _{j}\left( \omega
\right) {\bf 1}_{\left[ t_{j},t_{j+1}\right) }(t),
\end{equation}
where $\pi _{T}=\{0=t_{0}, t_1, \ldots ,t_{N}=T\}$ is a given partition of $\left[ 0,T%
\right] $ and $\xi _{i}\in Lip(\Omega
_{t_{i}})$, for all $i=0,1,2,\ldots ,N-1$. For $p\geq 1$ and $\eta \in \mathcal H_{G}^{0}\left( 0,T\right) $, define $%
\left\Vert \eta \right\Vert _{\mathcal{H}_{G}^{p}}=\{\mathbb{E}[(\int_{0}^{T}|\eta
_{s}|^{2}ds)^{p/2}]\}^{1/p}$, 
 and denote by $\mathcal H_{G}^{p}\left(
0,T\right) $ the completion of $%
\mathcal H_{G}^{0}\left( 0,T\right) $ under the norm $\left\Vert \cdot \right\Vert
_{\mathcal H_{G}^{p}}$.
\end{definition}
\begin{definition}[$G$-stochastic integrals]
For $\eta \in \mathcal{H}_{G}^{0}\left( 0,T\right) $ of the form (\ref{simp}),
the It\^{o} integral with respect to $G$-Brownian motion is defined by the linear mapping $%
\mathcal I: \mathcal{H}_{G}^{0}(0,T)\rightarrow L_{G}^{2}(\Omega _{T})$
\[
\mathcal I(\eta ):=\int_{0}^{T}\eta _{t}dB_{t}=\sum_{j=0}^{N-1}\xi
_{j}(B_{t_{j+1}}-B_{t_{j}}),\]%
which can be continuously extended to $I:\mathcal H_{G}^{1}(0,T)\rightarrow
L_{G}^{2}(\Omega _{T}).$ On the other hand, the stochastic integral with respect to $(\langle B\rangle_t)_{t\geq 0}$ is defined by the linear mapping $%
\mathcal Q: \mathcal{H}_{G}^{0}(0,T)\rightarrow L_{G}^{1}(\Omega _{T})$
\[
\mathcal Q(\eta ):=\int_{0}^{T}\eta _{t}d\left\langle B\right\rangle
_{t}=\sum_{j=0}^{N-1}\xi _{j}(\left\langle B\right\rangle
_{t_{j+1}}-\left\langle B\right\rangle _{t_{j}}),
\]%
which can also be continuously extended to $\mathcal{Q}: \mathcal H_{G}^{1}(0,T)\rightarrow
L_{G}^{1}(\Omega _{T}).$
\end{definition}

We have moreover some properties of the $G$-It\^o type integrals.
\begin{proposition}[B-D-G type inequality]
For $\eta \in \mathcal{H}^\alpha_G(0, T)$, $\alpha\geq 1$ and $p>0$, we have, 
\[
\underline{\sigma }^{p}c_{p}\mathbb{E}\left[ \left( \int_{0}^{T}|\eta
_{s}|^{2}ds\right) ^{p/2}\right] \leq \mathbb{E}\left[ \sup_{t\in \left[ 0,T%
\right] }\left\vert \int_{0}^{t}\eta _{s}dB_{s}\right\vert ^{p}\right] \leq
\overline{\sigma }^{p}C_{p}\mathbb{E}\left[ \left( \int_{0}^{T}|\eta
_{s}|^{2}ds\right) ^{p/2}\right] ,
\]%
where $0<c_{p}<C_{p}<\infty $ are constants independent of $\eta$,
$\underline{\sigma }$ and $\overline{\sigma }$.
\end{proposition}
\begin{proposition}
For all $\eta $, $\theta \in \mathcal{H}^{\alpha}_G\left( 0,T\right)$, $\alpha\geq 1$, with a bounded random variable $\xi\in L^1_G(\Omega_t)$, we
have
\[
\mathbb{E}_t\left[ \int_{t}^{T}\eta _{s}dB_{s}\right] =0;
\]%
\[
\int_{t}^{T}(\xi \eta _{s}+\theta _{s})dB_{s}=\xi \int_{t}^{T}\eta
_{s}dB_{s}+\int_{0}^{T}\theta _{s}dB_{s}.
\]%
\end{proposition}

Finally, we define the space $\mathcal{S}^p_G(0, T)$ for solutions of quadratic GBSDEs. Let
\begin{align*}
\mathcal S_{G}^{0}\left( 0,T\right) :=\big\{ h(t,B_{t_{1}\wedge t}, B_{t_{2}\wedge t}- B_{t_{1}\wedge t},&\ldots
,B_{t_{n}\wedge t} - B_{t_{n-1}\wedge t}):\\
&\ 0\leq t_{1}\leq t_2 \ldots ,t_{n}\leq T ,~h\in
\mathcal C_{b,Lip}(\mathbb{R}^{n+1})\big\},
\end{align*}
where $\mathcal C_{b,Lip}(\mathbb{R}^{n+1})$ is the collection of all bounded and Lipschitz functions on $\mathbb{R}^{n+1}$.
For $p\geq 1$ and $\eta \in \mathcal S_{G}^{0}\left( 0,T\right) $, we set $\left\Vert
\eta \right\Vert _{\mathcal S_{G}^{p}}=\{\mathbb{E}[\sup_{t\in \left[ 0,T\right]
}|\eta _{t}|^{p}]\}^{1/p}$. We denote by $\mathcal{S}_{G}^{p}\left( 0,T\right) $ the
completion of $\mathcal S_{G}^{0}\left( 0,T\right) $ under the norm $\left\Vert
\cdot\right\Vert _{\mathcal S_{G}^{p}}.$

%
%
%
%
%
%
%
%
%
%
\subsection{The formulation of GBSDEs}
In this paper, we shall consider the following type of equation:
\begin{equation}
Y_{t}=\xi +\int_{t}^{T}h(s, \omega_{\cdot\wedge s}, Y_{s},Z_{s})d\langle B\rangle
_{s}-\int_{t}^{T}Z_{s}dB_{s}-K_{T}+K_{t},  \label{qgbsde}
\end{equation}
where the terminal value $\xi$ and the generator $h$ satisfies the following conditions.
\begin{assumption}\label{assbig}
Assume that the generator $h: [0, T]\times\Omega \times \mathbb{R}^{2}\rightarrow \mathbb{R}$ satisfies the following conditions:\\
\noindent (\textbf{H0}) For each $(t,\omega )\in [0, T]\times \Omega$, $|h(t,\omega ,0,0)| +|\xi(\omega)|\leq M_0$;\\
\noindent (\textbf{Hc}) Moreover, $h(\cdot, \cdot, y, z)$ is uniformly continuous in $(t, \omega)$ and the modulus of continuity is independent of $(y, z)$, i.e., for each $(y, z)\times\mathbb{R}^2$,
$$\left|h(t_1 ,\omega^1 ,y,z)-h(t_2 ,\omega^2 ,y,z)\right|\leq {w}^h\left(|t_1-t_2|+\Vert\omega_1-\omega_2\Vert_\infty\right);$$
\noindent (\textbf{Hq}) The function $h$ is uniformly Lipschitz in $y$ and
uniformly locally Lipschitz in $z$, i.e., for each $(t,\omega )\in [0, T]\times \Omega$,
\begin{equation*}
|h(t,\omega ,y^{1},z^{1})-h(t,\omega ,y^{2},z^{2})|\leq
L_{y}|y^{1}-y^{2}|+L_{z}(1+|z^{1}|+|z^{2}|)|z^{1}-z^{2}|.
\end{equation*}
\end{assumption}
\begin{remark} From Theorem 4.7 in \cite{HWZ14}, we have for each $(y, z)\times\mathbb{R}^2$, $h(\cdot, \cdot, y, z)\in \mathcal{H}^2_G(0, T)$. Due to the boundedness, in fact, for any $p\geq 2$, $h(\cdot, \cdot, y, z)\in \mathcal{H}^p_G(0, T)$. Therefore, if $Y\in \mathcal{S}^p_G(0, T)$, $Z\in \mathcal{H}^{2p}_G(0, T)$,  $h(\cdot, \cdot, Y, Z)\in \mathcal{H}^p_G(0, T)$.
\end{remark}
\begin{remark}
One can always find a concave and sub-additive modulus $w$ in (Hc).
\end{remark}
\begin{definition}
For $p\geq 2$, a triple of processes $({Y}, {Z}, {K})$ belongs to $\mathfrak{G}^p_G(0, T)$, if $Y\in \mathcal{S}^p_G(0, T)$, $Z\in \mathcal{H}^p_G(0, T)$ and $K$ is a decreasing $G$-martingale with $K_0=0$ and $K_T\in L^p_G(\Omega_T)$. The triple $({Y}, {Z}, {K})$ is said solution of GBSDE (\ref{qgbsde}), if $({Y}, {Z}, {K})\in \mathfrak{G}^p_G(0, T)$, and for $0\leq t\leq T$, it satisfies (\ref{qgbsde}).
\end{definition}
\begin{remark}
All results in this paper hold for the GBSDE in a more general form, by trivially generalizing the argument:
$$
Y_{t}=\xi  +\int_{t}^{T}g(s, \omega_{\cdot\wedge s}, Y_{s},Z_{s})d{s}+\int_{t}^{T}h(s, \omega_{\cdot\wedge s}, Y_{s},Z_{s})d\langle B\rangle
_{s}-\int_{t}^{T}Z_{s}dB_{s}-K_{T}+K_{t},
$$
where the generator $g$ also satisfies Assumption \ref{assbig}.
\end{remark}
To prove the existence of solutions to (\ref{qgbsde}), we shall study the following auxiliary GBSDEs, in which the generator $h$ is replaced by a discrete function $f$, and correspondingly, a discrete terminal value is introduced here. Fixing $N\in \mathbb{N}$ and a partition $\pi^N:=\{0=t_0, t_1, \ldots, t_N=T\}$ on $[0, T]$, we consider
\begin{align}\label{po}
{Y}_t=\varphi(B_{t_1}, B_{t_2}&-B_{t_1}, \ldots, B_{t_{N}}-B_{t_{N-1}})\notag\\
&+\int^T_t {f}(s, B_{t_1\wedge s}, B_{t_2\wedge s}-B_{t_1\wedge s}, \ldots, B_{t_{N}\wedge s}-B_{t_{{N-1}}\wedge s}, {Y}_s, {Z}_s)d\langle B\rangle_s\\
&-\int^T_t {Z}_sdB_s-{K}_T+{K}_t.\notag
\end{align}
Similarly to Assumption \ref{assbig}, we introduce the following conditions on the terminal value $\varphi(B_{t_1}, B_{t_2}-B_{t_1}, \ldots, B_{t_N}-B_{t_{N-1}})$ and the generator $f$.
\begin{assumption}\label{apde}
We assume that the generator $f:[0, T]\times\mathbb{R}^N\times\mathbb{R}^2$ satisfies the following conditions:\\
\noindent (\textbf{H0'}) For each $(t, x_1, x_2, \ldots, x_N)\in [0, T]\times \mathbb{R}^N$, $|f(t,x_1, x_2,\ldots x_N ,0,0)|+|\varphi(x_1, x_2,\ldots x_N)| \leq M_0$;\\
\noindent (\textbf{Hc'}) Moreover, $f(\cdot, \cdot,\ldots, \cdot, y, z)$ is uniformly continuous in $(t, x_1, x_2, \ldots, x_N)$ and the modulus of continuity is independent of $(y, z)$, i.e., for each $(y, z)\in\mathbb{R}^2$,
$$\left|f(t_1 ,x^1_1, x^1_2, \ldots x^1_N, y, z)-f(t_1 ,x^2_1, x^2_2, \ldots x^2_N,y,z)\right|\leq {w}^f\left(|t_1-t_2|+\sum_{i=1}^N|x^1_i-x^2_i|\right);$$
\noindent (\textbf{Hq'}) The function $f$ is uniformly Lipschitz in $y$ and
uniformly locally Lipschitz in $z$, i.e., for each $(t, x_1, x_2, \ldots, x_N)\in [0, T]\times \mathbb{R}^N$,
\begin{equation*}
|f(t, x_1, x_2, \ldots x_N, y^{1},z^{1})-f(t,  x_1, x_2, \ldots x_N, y^{2},z^{2})|\leq
L_{y}|y^{1}-y^{2}|+L_{z}(1+|z^{1}|+|z^{2}|)|z^{1}-z^{2}|.
\end{equation*}
\end{assumption}

%

\section{$G$-Girsanov theorem and estimates for GBSDEs}
In this section, we first exhibit  the Girsanov type theorem in the $G$-framework by introducing the notion of $G$-BMO martingale generators. And then, we deduce a priori estimates for solutions of GBSDEs under quadratic assumptions.
\subsection{The Girsanov type theorem}
Similarly to Possama\"i and Zhou \cite{PZ13}, we can generalize the definition of BMO-martingale generators to the $G$-framework:
\begin{definition}
Assuming $Z\in \mathcal{H}^2_G(0. T)$, we say that $Z$ is a $G$-BMO martingale generator if
\begin{equation*}
\Vert Z\Vert _{BMO_G}^{2}:=\sup_{\mathbb{P}\in
\mathcal{P}_G}\left[ \sup_{\mathbb{\tau }\in \mathcal{T}_{0}^{T}}\left|
{E}_{\tau }^{\mathbb{P}}\left[ \int_{\tau }^{T}|
Z_{t}| ^{2}d\left\langle B\right\rangle _{t}\right] \right|
_{\infty }\right] <+\infty ,
\end{equation*}%
where $\mathcal{T}_{0}^{T}$ denotes the collection of all $\mathcal{F}$-stopping times taking values in $[0, T]$.
\end{definition}

\begin{lemma}\label{gsm}
Suppose that $Z\in \mathcal{H}_{G}^{2}( 0,T)$ is a G-BMO
martingale generator, then $$\mathcal{E}\left(\int Z dB\right):=\exp\left(\int Z dB-\frac{1}{2}\int |Z|^2 d\langle B\rangle\right)$$ is a symmetric $G$-martingale.
\end{lemma}

\noindent\textbf{Proof}: First, we verify that $\mathcal{E}_t(\int Z dB)\in L^1_G(\Omega_t)$ by Theorem 54 in Denis et al. \cite{DHP11}. As mentioned in Lemma 2.1 in \cite{PZ13}, if for some $q>1$ such that $\Vert Z\Vert_{BMO_G}\leq \Phi(q)$ (see Theorem 3.1 in Kazamaki \cite{Kaz94}), then
\begin{equation*}
{\mathbb{E}}\left[ \bigg( \mathcal{E}\bigg( \int_{0}^{\cdot
}Z_{s}dB_{s}\bigg)_t \bigg) ^{q }\right] <+\infty .
\end{equation*}
Fixing $t\in \left[ 0,T\right] $, we have for all $N>0$,
${\mathbb{E}}\left[ \mathcal{E}( Z)_t \mathbf{1}_{\left\{
\mathcal{E}(Z)_t \geq N\right\} }\right] \leq \frac{1}{%
N^{q-1 }}{\mathbb{E}}\left[ ( \mathcal{E}(Z)_t) ^{q}\right] $, which implies
$$\lim_{N\rightarrow \infty }{\mathbb{E}}\left[ \mathcal{E}( Z)_t 1_{\left\{ \mathcal{E}( Z)_t \geq N\right\}
}\right] =0.$$
Moreover, the quasi-continuity of $\mathcal{E}(\int Z dB)_t$ is inherited from $\int^t_0 Z_s dB_s$ and $\int^t_0 |Z|^2 d\langle B\rangle_s$. Thus, $\mathcal{E}(\int Z dB)_t\in L^1_G(\Omega_t)$.\\

From the well known results in \cite{Kaz94}, for $G$-BMO martingale generator $Z$, the process $\mathcal{E}(\int Z dB)$ is a martingale under each $\mathbb{P}\in\mathcal{P}_1$. 
By the representation of $G$-conditional expectation (Proposition 3.4 in Soner et al. \cite{STZ11}), we can deduce the desired result. \hfill{}$\square$\\

Following the procedure introduced in \cite{XSZ11}, we can define a new $G$-expectation on the space $Lip(\Omega_T)$ with $\mathcal{E}(Z)$ by
\begin{equation}\label{change}
\tilde{{\mathbb{E}}}[X]:=\sup_{\mathbb{P}\in\mathcal{P}_G}E^\mathbb{P}[\mathcal{E}(Z)_T X].
\end{equation}
Then, complete $Lip(\Omega_T)$ under $\tilde{{\mathbb{E}}}[\cdot]$ and obtain $\tilde{L}^{1}_{\tilde{G}}(\Omega_T)$.
If $X\in L^p_G(\Omega_T)$ with $p> \frac{q}{q-1}$, $ X\in \tilde{L}^{1}_{\tilde{G}}(\Omega_T)$. \\

The conditional expectation $\tilde{{\mathbb{E}}}_t[\cdot]$ thus can be first defined on $Lip(\Omega_T)$ then on $\tilde{L}^{1}_{\tilde{G}}(\Omega_T)$. 
 Obviously, for $G$-BMO martingale generator $Z$ and any $p\geq 1$, $\Vert Z\Vert_{\mathcal{H}^p_G(0, T)}<\infty$, and thus
$$
\left(\int^t_0Z_s d\langle B\rangle_s\right)^p\in L^1(\Omega_t)\ \ {\rm and} \ \ \left(\int^t_0Z_s d\langle B\rangle_s\right)^p\in \tilde{L}^1(\Omega_t).
$$

Notice that Xu et al. \cite{XSZ11} assumed the reinforced Novikov condition on $Z$ (Assumption 2.1 in \cite{XSZ11}) and develop the $G$-Girsanov type theory. This condition is mostly used for the proof of Lemma 2.2 in \cite{XSZ11} (corresponding to Lemma \ref{gsm} in the present paper) and thus, substituting this condition by a BMO one will not alter  the theory in \cite{XSZ11}. In particular, we have
\begin{lemma}
Suppose that $Z$ is a $G$-BMO martingale generator. We define a new $G$-expectation $\tilde{{\mathbb{E}}}[\cdot]$ by $\mathcal{E}(Z)$. Then, the process $B-\int Zd\langle B\rangle$ is a $G$-Brownian motion under $\tilde{\mathbb{E}}[\cdot]$.
\end{lemma}
The following result shows that a decreasing $G$-martingale under $\mathbb{E}[\cdot]$ is still a decreasing $G$-martingale under $\tilde{\mathbb{E}}[\cdot]$, if it satisfies a sufficient integrability condition.
\begin{lemma}\label{demart}
Assume the same as in the above lemma. Suppose that $K$ is a decreasing $G$-martingale such that $
K_{0}=0$ and for some $p>\frac{q}{q-1}$, $K_{t}\in L_{G}^{p}(\Omega _{t})$, $0\leq t\leq T$, where $q$ is the order in the reverse h\"older inequality for $\mathcal{E}(Z)$. Then $K$ is a
decreasing $G$-martingale under the new $G$-expectation defined by (\ref{change}).
\end{lemma}

\noindent\textbf{Proof :} The integrability condition on $K$ ensures that $\tilde{{\mathbb{E}}}_s[K_t]$, $0\leq s\leq t\leq T$, is well defined in the space $\tilde{L}_G^1(\Omega_s)$.
To prove this lemma, it suffices to verify the martingale property. Indeed, we recall Proposition 3.4 in Soner et al. \cite{STZ11} and deduce 
for some $0<\alpha<1$ such that $p>\frac{\alpha q}{\alpha q-1}$, and for $\mathbb{P}\in \mathcal{P}_G$, $\mathbb{P}$-a.s.,
\begin{eqnarray*}
0&\geq& \tilde{{\mathbb{E}}}_{t}\left[ K_{T}-K_{t}\right] \\[9pt]
&=&{\mathbb{E}}_{t}\left[
\frac{\mathcal{E}(Z)_T }{\mathcal{E}(Z)_t }%
( K_{T}-K_{t}) \right]= \mathop{\mbox{esssup}^\mathbb{P}}\limits_{\mathbb{P}'\in \mathcal{P}_G(t, \mathbb{P})} {E}^{\mathbb{P}'}_t\left[
\frac{\mathcal{E}(Z)_T }{\mathcal{E}(Z)_t }%
( K_{T}-K_{t}) \right]
=-\mathop{\mbox{essinf}^\mathbb{P}}\limits_{\mathbb{P}'\in \mathcal{P}_G(t, \mathbb{P})} {E}^{\mathbb{P}'}_t\left[
\frac{\mathcal{E}(Z)_T }{\mathcal{E}(Z)_t }%
( K_{t}-K_{T}) \right]\\
&\geq& -\mathop{\mbox{essinf}^\mathbb{P}}\limits_{\mathbb{P}'\in \mathcal{P}_G(t, \mathbb{P})}\left({{E}}^{\mathbb{P}'}_{t}\left[\left(
\frac{\mathcal{E}(Z)_T}{\mathcal{E}(Z)_t}\right)^{q} \right]^{\frac{1}{q}}
{{E}}^{\mathbb P'}_{t}\left[
( K_{t}-K_{T})^{\frac{\alpha q}{\alpha q-1}} \right]^{\frac{\alpha q-1}{q}}{{E}}^{\mathbb{P}'}_{t}\left[
 K_{t}-K_{T} \right]^{1-\alpha}\right)\\
&\geq&\mathop{\mbox{esssup}^\mathbb{P}}\limits_{\mathbb{P}'\in \mathcal{P}_G(t, \mathbb{P})} {E}^{\mathbb{P}'}_t\left[\left(
\frac{\mathcal{E}(Z)_T}{\mathcal{E}(Z)_t}\right)^{q} \right]^{\frac{1}{q}} \mathop{\mbox{esssup}^\mathbb{P}}\limits_{\mathbb{P}'\in \mathcal{P}_G(t, \mathbb{P})}  {{E}}^{\mathbb P'}_{t}\left[
( K_{t}-K_{T})^{\frac{\alpha q}{\alpha q-1}} \right]^{\frac{\alpha q -1}{q}}  \left(-\mathop{\mbox{essinf}^\mathbb{P}}\limits_{\mathbb{P}'\in \mathcal{P}_G(t, \mathbb{P})}{{E}}^{\mathbb{P}'}_{t}\left[
K_{t}-K_{T}\right]^{1-\alpha} \right)\\
&=&{\mathbb{E}}_{t}\left[\left(
\frac{\mathcal{E}(Z)_T}{\mathcal{E}(Z)_t}\right)^{q} \right]^{\frac{1}{q}}
{\mathbb{E}}_{t}\left[
( K_{t}-K_{T})^{\frac{\alpha q}{\alpha q-1}} \right]^{\frac{\alpha q -1}{q}}{\mathbb{E}}_{t}\left[
 K_{T}-K_{t} \right]^{1-\alpha}=0,
\end{eqnarray*}
where $\mathcal{P}_G(t, \mathbb{P}):=\{\mathbb{P}': \mathbb{P}'\in \mathcal{P}_G,\ \mathbb{P}'|_{\mathcal{F}_t}=\mathbb{P}|_{\mathcal{F}_t}\}$.
We apply the result of Proposition \ref{emu} to end the proof.
\hfill{}$\square$\\
\subsection{A priori estimates for GBSDEs} Consider the solution tripe $(Y, Z, K)$ for either (\ref{qgbsde}) or (\ref{po}). Applying Lemma 3.1 and Theorem 3.2 in \cite{PZ13}, we could easily obtain a upper bound for $Y$ and the $G$-BMO norm of $Z$, i.e.,
\begin{align}\label{estyz}
\big|\big|\sup_{0\leq t\leq T}|Y_t|\big|\big|_{L^{\infty}_G}+\Vert Z\Vert_{BMO_G} \leq \hat{C}:= C(M_0, L_y, L_z),
\end{align}
which implies that $Z$ is a $G$-BMO martingale generator. Moreover, for $p\geq 1$,
\begin{align}\label{estk}
\mathbb{E}[|K_t|^p]\leq \tilde{C}_p:= C(p, M_0, L_y, L_z),\quad 0\leq t\leq T.
\end{align}

Now we establish the following stability results for both (\ref{qgbsde}) and (\ref{po}) (compare with Theorem 3.2 in \cite{PZ13}).
\begin{proposition}\label{esty}Consider two quadratic GBSDEs (\ref{qgbsde}) with parameter $(\xi^1, h^1)$ and $(\xi^2, h^2)$, where $(\xi^i, h^i)$ satisfies (H0) and (Hq) with the same constants $M_0$, $L_y$ and  $L_z$.
Suppose that $(Y^i, Z^i, K^i)\in \mathfrak{G}^p_G(0, T)$, $p\geq 2$, are solutions corresponding to these parameters.
Then, for all $0\leq t\leq T$,%
\begin{equation*}
| Y_{t}^{1}-Y_{t}^{2}| \leq C\left( \Vert \xi ^{1}-\xi ^{2}\Vert_{L^\infty_G}+\tilde{\mathbb{E}}_t\left[ \int_{t}^{T}| h^{1}\left(
s,Y_{s}^{2},Z_{s}^{2}\right) -h^{2}\left( s,Y_{s}^{2},Z_{s}^{2}\right)
| d\langle B\rangle_s\right] \right),  \label{2.7}
\end{equation*}
where $\tilde{\mathbb{E}}[\cdot]$ is the new $G$-expectation under the Girsanov transform induced by $
\mathcal{E}(-b^\varepsilon)$ and $b^\varepsilon$ is defined in (\ref{poe}).
\end{proposition}
\noindent \textbf{Proof:} Let $\hat{Y}:=Y^{1}-Y^{2}$, $\hat{Z}:=Z^{1}-Z^{2}$, $\hat{K}%
:=K^{1}-K^{2}$ and $\hat{\xi}:=\xi ^{1}-\xi ^{2}$. For $0\leq t\leq T$, we have
\begin{align*}
\hat{Y}_{t} =\hat{\xi}&+\int_{t}^{T}\left(h^{1}\left(
s,Y_{s}^{1},Z_{s}^{1}\right) -h^{2}\left( s,Y_{s}^{2},Z_{s}^{2}\right) %
\right) d\left\langle B\right\rangle _{s}\\
&-\int_{t}^{T}\hat{Z}_{s}dB_{s}
-\left( K_{T}^{1}-K_{t}^{1}\right) +\left( K_{T}^{2}-K_{t}^{2}\right) .
\end{align*}%
We employ a linearization argument similar to the proof of Theorem 3.6 in \cite{HJPS14b}
by setting for $0\leq s\leq T$,
\begin{align}\label{poe}
\hat{a}_{s}^{\varepsilon } &=( 1-l( \hat{Y}_{s}) )
\frac{h^{1}\left( s,Y_{s}^{1},Z_{s}^{1}\right) -h^{1}\left(
s,Y_{s}^{2},Z_{s}^{1}\right) }{\hat{Y}_{s}}\mathbf{1}_{\{| \hat{Y}%
_{s}| >0\}}; \notag\\
\hat{b}_{s}^{\varepsilon } &=( 1-l(  \hat{Z}%
_{s} ) ) \frac{h^{1}\left(
s,Y_{s}^{2},Z_{s}^{1}\right) -h^{1}\left( s,Y_{s}^{2},Z_{s}^{2}\right) }{%
\vert \hat{Z}_{s}\vert ^{2}}\hat{Z}_{s}\mathbf{1}_{\{|
\hat{Z}_{s}| >0\}}; \\
\hat{m}_{s}^{\varepsilon } &=l( \hat{Y}_{s}) \left(h^{1}\left(
s,Y_{s}^{1},Z_{s}^{1}\right) -h^{1}\left( s,Y_{s}^{2},Z_{s}^{1}\right) %
\right)
+l( \hat{Z}_{s}) \left( h^{1}\left(
s,Y_{s}^{2},Z_{s}^{1}\right) -h^{1}\left( s,Y_{s}^{2},Z_{s}^{2}\right) %
\right); \notag\\
\hat{h}_{s} &=h^{1}\left( s,Y_{s}^{2},Z_{s}^{2}\right) -h^{2}\left(
s,Y_{s}^{2},Z_{s}^{2}\right),\notag
\end{align}
where $l$ is a Lipschitz function such
that $\mathbf{1}_{\left[ -\varepsilon ,\varepsilon \right] }( x)
\leq l( x) \leq \mathbf{1}_{\left[ -2\varepsilon ,2\varepsilon %
\right] }( x)$.
So we have%
\begin{equation*}
\hat{Y}_{t}=\hat{\xi}+\int_{t}^{T}\left( \hat{h}_{s}+\hat{m}%
_{s}^{\varepsilon }+\hat{a}_{s}^{\varepsilon }\hat{Y}_{s}+\hat{b}%
_{s}^{\varepsilon }\hat{Z}_{s}\right) d\left\langle B\right\rangle
_{s}-\int_{t}^{T}\hat{Z}_{s}dB_s-\int_{t}^{T}dK_{s}^{1}+%
\int_{t}^{T}dK_{s}^{2},
\end{equation*}
and
\begin{equation*}
| \hat{a}_{s}^{\varepsilon }| \leq L_{y},%
| \hat{b}_{s}^{\varepsilon }| \leq L_{z}\left(
1+| Z_{s}^{1}| +| Z_{s}^{2}| \right) ,%
\end{equation*}
\begin{equation*}| \hat{m}_{s}^{\varepsilon }| \leq
2\varepsilon \left(L_{y}+L_{z}\left( 1+2\varepsilon +2|
Z_{s}^{1}| \right) \right).
\end{equation*}
From the proof of Theorem 3.6 in Hu et al. \cite{HJPS14b},
we know the process $\hat{a}^{\varepsilon }$
belongs to $\mathcal{H}^2_G(0, T)$ and furthermore due to the boundedness, $\hat{a}^\varepsilon$ and $e^{\int \hat a^\varepsilon d\langle B\rangle}$ belong to $\mathcal{H}^p_G(0, T)$, for any $p\geq 2$.\\

Applying It\^o's formula to $e^{\int \hat a^\varepsilon d\langle B\rangle}\hat{Y}$, we have
\begin{align*}
e^{\int^t_0 \hat a^\varepsilon_s d\langle B\rangle_s}\hat{Y}_{t}=e^{\int^T_0 \hat a^\varepsilon_t d\langle B\rangle_t}\hat{\xi}&+\int_{t}^{T}e^{\int^s_0 \hat a^\varepsilon_u d\langle B\rangle_u}\left( \hat{h}_{s}+\hat{m}%
_{s}^{\varepsilon }+\hat{b}%
_{s}^{\varepsilon }\hat{Z}_{s}\right) d\left\langle B\right\rangle
_{s}\\&-\int_{t}^{T}e^{\int^s_0 \hat a^\varepsilon_u d\langle B\rangle_u}\hat{Z}_{s}dB_s-\int_{t}^{T}e^{\int^s_0 \hat a^\varepsilon_u d\langle B\rangle_u}dK_{s}^{1}+%
\int_{t}^{T}e^{\int^s_0 \hat a^\varepsilon_u d\langle B\rangle_u}dK_{s}^{2}.
\end{align*}
Furthermore,
\begin{align*}
e^{\int^t_0 \hat a^\varepsilon_s d\langle B\rangle_s}\hat{Y}_{t}+\int_{t}^{T}e^{\int^s_0 \hat a^\varepsilon_u d\langle B\rangle_u}dK_{s}^{1}\leq e^{\int^T_0 \hat a^\varepsilon_t d\langle B\rangle_t}\hat{\xi}
&+\int_{t}^{T}e^{\int^s_0 \hat a^\varepsilon_u d\langle B\rangle_u}\left( \hat{h}_{s}+\hat{m}%
_{s}^{\varepsilon }+\hat{b}%
_{s}^{\varepsilon }\hat{Z}_{s}\right) d\left\langle B\right\rangle
_{s}\\&-\int_{t}^{T}e^{\int^s_0 \hat a^\varepsilon_u d\langle B\rangle_u}\hat{Z}_{s}dB_s.
\end{align*}

Thanks to the lemma below, $\hat{b}^\varepsilon$ belongs to $\mathcal{H}^2_G(0, T)$. Moreover, $\hat{b}^\varepsilon$ is a $G$-BMO martingale generator. Thus, we could define a new $G$-expectation $\hat{\mathbb{E}}[\cdot]$ by $\mathcal{E}(\hat b^\varepsilon)$, such that $\hat{B}:= B-\int \hat b^\varepsilon d\langle B\rangle$ is a $G$-Brownian motion under $\hat{\mathbb{E}}[\cdot]$. From (\ref{estk}), we know that $K_t$ is $p$-integrable, for any $p\geq 1$, thus we could assume without loss of generality that $p>\frac{q}{q-1}$ and $q$ satisfies $L_z(1+2\hat{C})<\Phi(q)$.  Since $e^{\int^\cdot_0 \hat a^\varepsilon_t d\langle B\rangle_t}$ is a positive process, from Lemma 3.4 in Hu et al. \cite{HJPS14a}, we know $\int^{\cdot}_{0}e^{\int^s_0 \hat a^\varepsilon_u d\langle B\rangle_u}dK_{s}^{1}$ and $\int_{0}^{\cdot}e^{\int^s_0 \hat a^\varepsilon_u d\langle B\rangle_u}dK_{s}^{2}$ is a decreasing $G$-martingale under both $\mathbb{E}[\cdot]$ and $\hat{\mathbb{E}}[\cdot]$.\\

Taking the conditional $G$-expectation $\hat{\mathbb{E}}_t[\cdot]$ on both sides, it is easy to see that
\begin{align*}
e^{\int^t_0 \hat a^\varepsilon_s d\langle B\rangle_s}\hat{Y}_{t}
&\leq \hat{\mathbb{E}}_t\bigg[
e^{\int^T_0 \hat a^\varepsilon_t d\langle B\rangle_t}\hat{\xi}
+\int_{t}^{T}e^{\int^s_0 \hat a^\varepsilon_ud\langle B\rangle_u}\left( \hat{h}_{s}+\hat{m}%
_{s}^{\varepsilon }\right) d\left\langle B\right\rangle
_{s}\bigg]\\
&\leq e^{{\overline{\sigma}}^2L_yT}\Vert\hat{\xi}\Vert_{L^\infty_G}+e^{{\overline{\sigma}}^2L_yT}\left(\hat{\mathbb{E}}_t\left[\int_{t}^{T}|\hat{h}_{s}|d\langle B\rangle_s\right]\right)\\
&+2\varepsilon e^{{\overline{\sigma}}^2L_yT}\left(\hat{\mathbb{E}}_t\left[\int_{t}^{T}\left(L_{y}+L_{z}\left(1+2\varepsilon +2|Z_{s}^{1}| \right)\right)d\langle B\rangle_s\right]\right),
\end{align*}
which implies
\begin{align*}
\hat{Y}_{t}
\leq e^{2{\overline{\sigma}}^2L_yT}\Vert\hat{\xi}\Vert_{L^\infty_G}&+e^{2{\overline{\sigma}}^2L_yT}\left(\hat{\mathbb{E}}_t\left[\int_{t}^{T}|\hat{h}_{s}|d\langle B\rangle_s\right]\right)\\
&+2\varepsilon e^{2{\overline{\sigma}}^2L_yT}\left(\hat{\mathbb{E}}_t\left[\int_{t}^{T}\left(L_{y}+L_{z}\left(1+2\varepsilon +2|Z_{s}^{1}| \right)\right)d\langle B\rangle_s\right]\right).
\end{align*}
Sending $\varepsilon \rightarrow 0$, we have
\begin{align*}
\hat{Y}_{t}
\leq e^{2{\overline{\sigma}}^2L_yT}\Vert\hat{\xi}\Vert_{L^\infty_G}&+e^{2{\overline{\sigma}}^2L_yT}\left(\hat{\mathbb{E}}_t\left[\int_{t}^{T}|\hat{h}_{s}|d\langle B\rangle_s\right]\right).\end{align*}
%
Moreover, we deduce in the same way that
\begin{align*}
-\hat{Y}_{t}
\leq e^{2{\overline{\sigma}}^2L_yT}\Vert\hat{\xi}\Vert_{L^\infty_G}&+e^{2{\overline{\sigma}}^2L_yT}\left(\hat{\mathbb{E}}_t\left[\int_{t}^{T}|\hat{h}_{s}|d\langle B\rangle_s\right]\right).\end{align*}
This completes the proof.
\hfill$\square$
\begin{lemma}
The process defined in (\ref{poe}) belongs to $\mathcal{H}^2_G(0, T)$ and is a $G$-BMO martingale generator.
\end{lemma}
\noindent {\bf Proof:} For each $n\in \mathbb{N}$, define $h^i_{n}$ as follows
\begin{equation*}
h^1_{n}( s,y,z) =h^1\left( s,y,\frac{| z| \wedge n%
}{| z| }z\right),
\end{equation*}%
which is Lipschitz in $z$ with the Lipschitz constant $L_{z}(1+2n)$. For
each $n\in \mathbb{N}$, the process $\hat{b}^{\varepsilon ,n}$ is defined by
\begin{equation*}
\hat{b}_{s}^{\varepsilon ,n}=( 1-l(\hat{Z}_{s} )) \frac{h^1_{n}( s,Y^{2}_s,Z^{1}_{s}) -h^1_{n}( s,Y^{2}_s,Z^{2}_s) }{\vert \hat{Z}_{s}\vert ^{2}}{%
\hat{Z}_{s}}\mathbf{1}_{\left\{ |\hat{Z}_{s}|>
0\right\} },\quad 0\leq s\leq T,
\end{equation*}%
which belongs to $\mathcal{H}^2_G(0, T)$ according to \cite{HJPS14b}.
On the other hand,
\begin{align*}
| \hat{b}_{s}^{\varepsilon ,n}-\hat{b}_{s}^{\varepsilon }| &\leq(
1-l(  \hat{Z}_{s} ) ) \frac{|
h^1( s,Y^2_s, Z^1_{s}) -h^1_{n}( s,Y^2_s,Z^1_{s}) | + |
h^1( s,Y^2_s, Z^2_{s}) -h^1_{n}( s,Y^2_s,Z^2_{s}) |}{%
 |\hat{Z}_{s}| }\mathbf{1}_{\left\{ |\hat{Z}_{s}|> 0\right\} } \\
&\leq \frac{L_z}{\varepsilon}\left((|Z^1_{s}|-n)(1+n+|Z^1_s|)\mathbf{1}_{\{ |Z^1 _{s}|>n\} }+(|Z^2_{s}|-n)(1+n+|Z^2_s|)\mathbf{1}_{\{ | Z^2_{s}|>n\} }\right)\\
&\leq \frac{CL_z}{\varepsilon}\left (|Z^1_{s}|^2\mathbf{1}_{\left\{ | Z^1_{s}|>n\right\} }+|Z^2_{s}|^2\mathbf{1}_{\{ | Z^2_{s}|>n\} }\right),
\end{align*}%
where $C$ is independent of $n$. By Proposition 2.9 in Li and Peng \cite{LP11},  we conclude that $\hat{b}^\varepsilon$ belongs to $\mathcal{H}^2_G(0, T)$. From the estimate for $\hat{b}^{\varepsilon}$, we deduce that  $\Vert\hat{b}^{\varepsilon}\Vert_{BMO_G}\leq L_z(1+2\hat{C})$. Thus, $\hat{b}^{\varepsilon}$ is a $G$-BMO
martingale generator.
\hfill$\square$

\begin{proposition}\label{estz}
Consider two quadratic GBSDEs (\ref{qgbsde}) with parameter $(\xi^1, h^1)$ and $(\xi^2, h^2)$, where $(\xi^i, h^i)$ satisfies (H0) and (Hq) with the same constants $M_0$, $L_y$ and  $L_z$.
Suppose $(Y^i, Z^i, K^i)\in \mathfrak{G}^p_G(0, T)$, $p\geq 2$, are solutions corresponding to these parameters.
Then, for $1\leq p'/2<p$,
$$
\mathbb{E}\left[\left(\int^T_0|Z^1_t-Z^2_t|^2dt\right)^{p'/2}\right]
\leq C({p, \overline{\sigma}, \underline{\sigma}}, M_0, L_y, L_z)\bigg(\Vert\xi^1-\xi^2\Vert^{p'}_{L^\infty_G}+\mathbb{E}\left[\sup_{0\leq t\leq T}|Y^1_t-Y^2_t|^p\right]^{p'/2p}\bigg).
$$
\end{proposition}
\noindent {\bf Proof:} We keep the notations in the proof of Proposition \ref{esty}. Indeed, due to the boundedness of the $G$-BMO norm of $Z^1$ and $Z^2$, $Z^1$, $Z^2\in \mathcal{H}^p_G(0, T)$, for any $p\geq 2$. Applying It\^o's formula to $|\hat Y|^2$, we have
\begin{align*}
2|\hat Y_0|^2 +\int^T_0|\hat Z_t|^2d\langle B\rangle_t\leq 2|\hat \xi|^2&+4 \int^T_0\hat Y_t \left(F^1(t, Y^1_t, Z^1_t)-F^2(t, Y^2_t, Z^2_t)\right)d\langle B\rangle_t\\
&-4\int^T_0\hat Y_t\hat Z_t dB_t-4\int^T_0\hat Y_{t}d \hat K_t.
\end{align*}
Then, for $1\leq p'/2<p$,
\begin{align*}
&\ \ \ \ \ \ \mathbb{E}\left[\left(\int^T_0|\hat Z_t|^2dt\right)^{p'/2}\right] \leq C({p, \overline{\sigma}, \underline{\sigma}})\bigg(\Vert\hat \xi\Vert^{p'}_{L^\infty_G}\\
&+\mathbb{E}\bigg[\bigg(\sup_{0\leq t\leq T}|\hat Y_t| \int^T_0 \left(2M_0+L_y\left(|Y^1_t|+|Y^2_t|\right)+2L_z\left(1+|Z^1_t|^2+|Z^2_t|^2\right)\right)d\langle B\rangle_t\bigg)^{p'/2}\bigg]\\
&+\mathbb{E}\bigg[\left(\sup_{0\leq t\leq T}|\hat Y_t|^2\int^T_0|\hat Z_t|^2d\langle B\rangle_t\right)^{p'/4}\bigg]+\mathbb{E}\bigg[\left(\sup_{0\leq t\leq T}|\hat Y_t||K^1_T|\right)^{p'/2}\bigg]\\
&+\mathbb{E}\bigg[\left(\sup_{0\leq t\leq T}|\hat Y_t||K^2_T|\right)^{p'/2}\bigg]\bigg).
\end{align*}
From (\ref{estyz}) and (\ref{estk}) and by H\"older's inequality, we have
\begin{align*}
\mathbb{E}\left[\left(\int^T_0|\hat Z_t|^2dt\right)^{p'/2}\right]&\leq C({p, \overline{\sigma}, \underline{\sigma}}, M_0, L_y, L_z)\bigg(\Vert\hat \xi\Vert^{p'}_{L^\infty_G}+\mathbb{E}\left[\sup_{0\leq t\leq T}|\hat{Y}_t|^p\right]^{p'/2p}\bigg).
\end{align*}
\hfill$\square$
\begin{remark}
The uniqueness for the quadratic GBSDE can be derived from Proposition \ref{esty} and Proposition \ref{estz}, or by regarding it as a quadratic 2BSDEs studied in \cite{PZ13}.
\end{remark}

\section{The existence of solutions to discrete GBSDEs}
In this section, we prove the existence of solutions to the equation (\ref{po}), which are constructed by solutions of the corresponding discrete PDEs.
\subsection{Discrete PDEs}
In this subsection, we follow Hu and Ma \cite{HM04} to consider discrete PDEs and deduce the boundedness of the first derivatives of the solution $u$ in ${x}$.

First, we introduce  the following fully nonlinear PDE on $[t_{N-1}, T]$:
\begin{equation}\label{pden}
\partial _{t}u+G(D_{x^Nx^N}^{2}u+2f(t,x_1, x_2, \ldots, x_N,u, D_{x_N}u))=0,
\end{equation}%
with $u(T,x_1, x_2, \ldots, x_N)=\varphi (x_1, x_2, \ldots, x_N)\in \mathcal{C}_{b, lip}(\mathbb{R}^N)$, where the terminal value $\varphi$ and the generator $f$ satisfies Assumption \ref{apde} and the following assumption.
\begin{assumption}\label{bpde} We assume that the generator $f:[0, T]\times\mathbb{R}^N\times\mathbb{R}^2$ satisfies moreover the following conditions:\\
\noindent (\textbf{Hd'}) The function $f$ is at least $\mathcal{C}^1$ in $x_1, x_2,\ldots, x_N, y, z$, differentiable in $t$ and twice differentiable in $x_1, x_2,\ldots, x_N, y, z$, where the second derivative of $f$ in $x_1, x_2,\ldots, x_N, y, z$ are bounded on the set $[0, T]\times \mathbb{R}^N\times  [-M_y, M_y]\times [-M_z, M_z]$, for any $M_y$, $M_z>0$.
\end{assumption}
\begin{remark}
From (Hc') and (Hq'), we could conclude that the first derivative of $f$ in $t$, $x_1, x_2,\ldots, x_N, y, z$ are bounded on the set $[0, T]\times \mathbb{R}^N\times \mathbb{R}\times [-M_z, M_z]$, for any $M_z>0$.
\end{remark}
\begin{proposition}
The PDE (\ref{pden}) admits a classical solution bounded by $M:=M(M_0, L_y)$, and for arbitrary small $\kappa$, there exists a constant $\alpha:=\alpha(\kappa)$, such that
\begin{equation*}\label{est}
\Vert u\Vert_{\mathcal{C}^{1+\alpha /2,2+\alpha }([t_{N-1},T-\kappa ]\times \mathbb{R})}<\infty.
\end{equation*}
\end{proposition}
\begin{remark} The proof of this proposition is not difficult by proceeding a similar argument as Appendix \S B-4 in Peng \cite{Pen10}, where the results from Example 6.1.8 and Theorem 6.4.3 in Krylov \cite{Kry87} play very important roles.
\end{remark}
Fix $N\in \mathbb{N}$ and a partition $\pi^N$ on $[0, T]$. Denote by $L^{\varphi}$ the Lipschitz constant of $\varphi $, and by $L^f_x
$ and $L^f_y$ the Lipschitz constants of $f$ in ${x_1}$,  ${x_2}\ldots {x_N}$ and in $y$.\\

Denote $\mathbf{x}^{(k)}:=(x_1, x_2, \ldots, x_k)$, $k=1, 2, \ldots, N$. We rewrite (\ref{pden}) into the following form:
\begin{align*}
\partial _{t} & u^{N} (t,\mathbf{x}^{(N-1)},x_N) +G(D_{x_Nx_N}^{2}u^{N}(t,\mathbf{x}%
^{(N-1)},x_N) \\
&+2f(t,\mathbf{x}^{(N-1)},x_N,u^{N}(t,\mathbf{x}^{(N-1)},x_N),D_{x_N}u^{N}(t,%
\mathbf{x}^{(N-1)},x_N)))=0,\notag
\end{align*}%
with $u^{N}(T,\mathbf{x}^{(N-1)},x)=\varphi (\mathbf{x}^{(N-1)},x)$.\\

To estimate the first derivative $\frac{\partial u}{\partial x_N}$, we proceed the same as Step 1 of proof for Theorem 4.1 in Hu et al. \cite{HJPS14a} and obtain
\begin{equation*}
\left\vert D_{x_N}u^{N}(t,\mathbf{x}^{(N-1)},x_N)\right\vert \leq \left(L^\varphi+%
\frac{L^f_x}{L^f_y}\right)\exp (\overline{\sigma }^{2}L_{y}^{f}\left(
T-t_{N-1}\right) )-\frac{L^f_x}{L^f_y}:=L^N.
\end{equation*}
\begin{remark} We remark here for proving the above result, we shall recall a general comparison result in Buckdhan and Li \cite{BL08}  (Theorem 6.1). Note that
the function $f$ we consider is only local Lipschitz, while Theorem 6.1 in \cite{BL08} requires the Lipschitz assumption. This has little matter, since we could eventually see that $D_{x_N}u$ is bounded by a constant $M_z:=M_z(N, L^\varphi, L_x, L_y)$, then a standard truncation technique may apply here.
\end{remark}
In a similar way, we have moreover, for $k=1, 2, \ldots, N$, $t\in [t_{N-1}, T]$,
\begin{equation*}
\left\vert D_{x_k}u^{N}(t,\mathbf{x}^{(N-1)}, x_N)\right\vert \leq L^N.
\end{equation*}

 Then, we could define the following PDE on $[t_{N-2},t_{N-1}]$:
\begin{align*} \label{e123}
\partial _{t} & u^{N-1} (t,\mathbf{x}^{(N-2)},x_{N-1}) +G(D_{x_{N-1}x_{N-1}}^{2}u^{N-1}(t,\mathbf{x}%
^{(N-2)},x_{N-1}) \\
&+2f(t,\mathbf{x}^{(N-2)},x_{N-1},0, u^{N-1}(t,\mathbf{x}^{(N-2)},x_{N-1}),D_{x}u^{N-1}(t,%
\mathbf{x}^{(N-2)},x_{N-1})))=0,\notag
\end{align*}%
with the terminal condition
\begin{equation*}
u^{N-1}(t_{N-1},\mathbf{x}^{(N-2)},x_{N-1}):=u^{N}(t_{N-1},\mathbf{x}^{(N-2)},x_{N-1},0).
\end{equation*}%
From the estimate above, we know that the Lipschitz constant of $u^N$ in $x_{N-1}$ is $L^N$, then for $k=1, 2, \ldots, N-1$, $t\in [t_{N-2}, t_{N-1}]$,
\begin{equation*}
\left\vert D_{x_k}u^{N-1}(t,\mathbf{x}^{(N-2)},x_{N-1})\right\vert \leq
\left(L^N+\frac{L^{f}_{x}}{L^{f}_{y}}\right)\exp (\overline{\sigma }%
^{2}L_{y}^{f}\left( t_{N-1}-t_{N-2}\right) )-\frac{L^{f}_{x}}{L^{f}_{y}}:=L^{N-1}.
\end{equation*}

By recurrence, we consider the following PDE on $[t_{k-1},t_{k}]$:
\begin{align}\label{pder}
\partial _{t}  u^{k} (t,\mathbf{x}^{(k-1)},x_{k}) &+G(D_{x_{k}x_{k}}^{2}u^{k}(t,\mathbf{x}%
^{(k-1)},x_{k}) \\
&+2f(t,\mathbf{x}^{(k-1)},x_{k},\underbrace{0, \ldots, 0}_{N-k}, u^{k}(t,\mathbf{x}^{(k-1)},x_{k}),D_{x_{k}}u^{k}(t,%
\mathbf{x}^{(k-1)},x_{k})))=0,\notag
\end{align}%
with the terminal condition
\begin{equation*}
u^{k}(t_{k},\mathbf{x}^{(k-1)},x_{k}):=u^{k+1}(t_{k},\mathbf{x}^{(k-1)},x_{k},0).
\end{equation*}%
For each $i=0, 1, 2, \ldots, k$, $t\in [t_{k-1}, t_{k}]$,
\begin{equation*}
\left\vert D_{{x}_{i}}u^{k}(t,\mathbf{x}^{(k-1)}, x_{k})\right\vert \leq
\left(L^{k+1}+\frac{L^{f}_{x}}{L^{f}_{y}}\right)\exp (\overline{\sigma }%
^{2}L_{y}^{f}\left( t_{k}-t_{k-1}\right))-\frac{L_{h}^{x}}{L_{h}^{y}}:=L^{k}.
\end{equation*}
\subsection{The solution of the discrete GBSDE}
In this subsection, we construct the solution of the discrete GBSDE (\ref{po}) satisfying Assumption \ref{apde} and \ref{bpde}. Fix $N\in \mathbb{N}$ and the partition $\pi^N$ of $[0, T]$. We note
$$\mathbf{B}_{t}^{k}:=\left(
B_{t_{1}},B_{t_{2}}-B_{t_{1}},\ldots
,B_{t_{k-1}}-B_{t_{k-2}},B_{t}-B_{t_{k-1}}\right),\ t\in[t_{k-1}, t_k],\ k=1, 2, \ldots, N.$$
Then, the solution of the GBSDE (\ref{po}) is defined in the following way: for $t\in[t_{k-1}, t_k]$,
\begin{equation*}
Y_{t}:=Y_{t}^{k}:=u^{k}(t,\mathbf{B}_{t}^{k});
\end{equation*}%
\begin{equation*}
Z_{t}:=Z_{t}^{k}:=D_{x_{k}}u^{k}(t,\mathbf{B}_{t}^{k});
\end{equation*}%
and
\begin{align*}
K_{t}:=K_{t}^{k}& =K_{t_{k-1}}^{k-1}+\int_{t_{k-1}}^{t}\left(
D_{x_{k}x_{k}}^{2}u^{k}(s,\mathbf{B}_{s}^{k})+2f(s,\mathbf{B}
_{s}^{k}, \underbrace{0, \ldots, 0}_{N-k}, u^{k}(s,\mathbf{B}_{s}^{k}),D_{x_{k}}u^{k}(s,\mathbf{B}%
_{s}^{k}))\right) d\langle B\rangle _{s} \\
& -\int_{t_{k-1}}^{t}G\left( D_{x_{k}x_{k}}^{2}u^{k}(s,\mathbf{B}%
_{s}^{k})+2f( s,\mathbf{B}_{s}^{k},\underbrace{0, \ldots, 0}_{N-k}, u^{k}(s,\mathbf{B}%
_{s}^{k}),D_{x_{k}}u^{k}(s,\mathbf{B}_{s}^{k}) \right) ds,
\end{align*}%
where $u^k$, $k=1, 2, \ldots, N$, is the solution to the corresponding PDE (\ref{pder}) on $[t_{k-1}, t_k]$.\\

From the estimates of the corresponding PDEs, we can find the following bounds:
$$
|Y|\leq M_y:=M_y(N, M_0, L^f_y)
$$
and
\begin{align}\label{solbou}
|Z|\leq M_z:=M_z(N, L^\varphi, L^f_x, L^f_y).
\end{align}

If $N=1$, then following the proof of Theorem 4.1 in \cite{HJPS14a}, we can
verify that for each $\kappa\in (0, 1)$, $(Y, Z, K)$ is a solution of the following
GBSDE on $[0, T-\kappa]$:
\begin{equation*}
Y_{t}=\varphi (B_{T-\kappa})+\int_{t}^{T-\kappa}f(t,B_{s},Y_{s},Z_{s})d\langle B\rangle
_{s}-\int_{t}^{T-\kappa}Z_{s}dB_{s}-(K_{T-\kappa}-K_{t}),
\end{equation*}%
where $K$ is a decreasing $G$-martingale on $[0, T-\kappa]$ with $K_{0}=0$. Indeed, the function $u$ and the derivative $Du$ and $D^2u$ is bounded and $\alpha$-H\"older continuous, and thus $D^2_{xx}u(\cdot, B)+2f(\cdot, B, u(\cdot, B), D_xu(\cdot, B))\in \mathcal{H}^p_G(0, T-\kappa)$, for any $p\geq 2$.\\

Similarly to (4.3) in \cite{HJPS14a}, we could obtain for $0<\tilde{t}\leq \hat{t}<T$ and $\tilde{x}$, $\tilde{x'}\in \mathbb{R}$ and some positive constant $L$,
$$|u(\tilde{t}, \tilde{x})-u(\hat{t},\hat{x})|\leq L\left(\sqrt{|\tilde{t}-\hat{t}|}+|\tilde{x}-\hat{x}|\right),$$
which implies that $(Y, Z, K)$ is the solution of (\ref{po}) on $[0, T]$, and $K$ is a decreasing $G$-martingale with $K_0=0$ and closed by $K_T\in L^p_G(\Omega_T)$, for any $p\geq 1$, which can be defined as the quasi-sure limit of the decreasing sequence $K_{T-\frac{1}{n}}$. We remark also that $Y\in \mathcal{S}^p_G(0, T)$ and $Z\in \mathcal{H}^p_G(0, T)$, for any $p\geq 2$, which could be deduced by the same procedure as \cite{HJPS14a}.

For the case that $N>1$, it suffices to prove without loss of generality that $(Y, Z, K)\in \mathfrak{G}^p_G(0, T)$ is the solution of (\ref{po}) on $[0, T]$ when $N=2$. \\

On $[0, t_1]$, it follows from the previous case $N=1$ that the triple $(Y, Z, K)$ defined by
$$
Y_\cdot := u^1(\cdot, B_\cdot),\quad Z_\cdot:= D_{x_1}u^1(\cdot, B_\cdot),$$
\begin{align*}
K_\cdot:=\int_{0}^{\cdot}&\left(
D_{x_{1}x_{1}}^{2}u^{1}(s, {B}_{s})+2f(s,B_{s}, 0, u^{1}(s,{B}_{s}),D_{x_{1}}u^{1}(s,{B}_{s}))\right) d\langle B\rangle_{s}\\
&-\int_{0}^{\cdot}G\left(
D_{x_{1}x_{1}}^{2}u^{1}(s, {B}_{s})+2f(s,B_{s}, 0, u^{1}(s,{B}_{s}),D_{x_{1}}u^{1}(s,{B}_{s}))\right)ds,
\end{align*}
solves the following GBSDE:
\begin{align*}
{Y}_t=u^2(t_1, B_{t_1}, 0)
&+\int^{t_1}_t {f}(s, B_{s}, 0, {Y}_s, {Z}_s)d\langle B\rangle_s-\int^{t_1}_t {Z}_sdB_s-({K}_{t_1}-{K}_t).
\end{align*}
%
%
%

Now, it suffice to verify that $$%
Y_{t}^{2}:=u^{2}\left( t,B_{t_{1}},B_{t}-B_{t_{1}}\right), \quad%
Z_{t}^{2}:=D_{x_{2}}u^{2}\left( t,B_{t_{1}},B_{t}-B_{t_{1}}\right),$$
\begin{align*}
K_{t}
:=K_{t_{1}}&+\int_{t_1}^{t}D_{x_{2}x_{2}}^{2}u^{2}(s,B_{t_{1}},B_{s}-B_{t_{1}}) \\
&
+2f(s,B_{t_{1}},B_{s}-B_{t_{1}}, u^{2}(s,B_{t_{1}},B_{s}-B_{t_{1}}),D_{x_{2}}u^{2}(s,B_{t_{1}},B_{s}-B_{t_{1}}))d\langle B\rangle _{s}
\\
& -\int_{t_{1}}^{t}G(D_{x_{2}x_{2}}^{2}u^{2}(s,B_{t_{1}},B_{s}-B_{t_{1}})
\\
&
+2f(s,B_{t_{1}},B_{s}-B_{t_{1}}, u^{2}(s,B_{t_{1}},B_{s}-B_{t_{1}}),D_{x_{2}}u^{2}(s,B_{t_{1}},B_{s}-B_{t_{1}})))ds,
\end{align*}%
defines a solution on $[t_{1},T]$ of the following GBSDE:
\begin{equation*}
Y_{t}=\varphi (B_{t_{1}},B_{T}-B_{t_{1}})+\int_{t}^{T}f(s,B_{t_{1}},B_{s\wedge t_{1}}-B_{t_{1}}, Y_{s},Z_{s})d\langle B\rangle
_{s}-\int_{t}^{T}Z_{s}dB_{s}-(K_{T}-K_{t}).
\end{equation*}%
This can be achieved by first proving a generalized It\^o's formula for $u(t, B_{t_1}, B_{t}-B_{t_1})$ on $[t_1, T-\kappa]$ as \S III-6 in \cite{Pen10} and then letting $\kappa\rightarrow 0$.
\section{Existence of solutions for general quadratic GBSDEs}
In this section, we shall prove the existence result for the general quadratic GBSDE (\ref{qgbsde}) under Assumption \ref{assbig}. Indeed, we start by considering the solution to the discrete GBSDE (\ref{po}) under assumptions weaker than (Hd') and construct solutions to (\ref{qgbsde}) by successive approximation.\\

\noindent {\bf Step 1}:
Fix $N\in \mathbb{N}$ and the partition $\pi^N$ on $[0, T]$. We consider the GBSDE (\ref{po}) with the generator $\hat{f}$ satisfying Assumption \ref{apde} and what follows.
\begin{assumption}\label{cpde} We assume that the generator $\hat f:[0, T]\times\mathbb{R}^N\times\mathbb{R}^2$ satisfies moreover the following conditions:\\
\noindent {\bf (Hd'')}\  The first derivative of $f$ in $t$, the first and the second derivatives of $f$ in $x_1, x_2,\ldots, x_N$ are bounded on the set $[0, T]\times \mathbb{R}^N\times [-{M}_y, {M}_y]\times [-{M}_z, {M}_z]$, for any ${M}_y$, ${M}_z>0$.
\end{assumption}

In what follows, we shall regularize $\hat{h}$ in $y$ and $z$: for each $(t, x_1, x_2,\ldots, x_N)\in [0, T]\times\mathbb{R}^{N}$, $(y, z)\in \mathbb{R}^2$, we define
\begin{align*}
{f}^{n}\left(t, x_1, x_2,\ \ldots, x_N,y,z \right) :=\int_{\mathbb{R}^2} \hat{f}\left( t, x_1, x_2, \ldots, x_N,y-\tilde{y},z-\tilde{z}\right) \rho_{n}\left(\tilde{y}, \tilde{z}\right) d\tilde{y}d\tilde{z},
\end{align*}
where $\rho_n$ is a positive smooth function such that its support is contained in a $\frac{1}{n}$-ball in $\mathbb{R}^{2}$ and $\int_{\mathbb{R}^{2}} \rho_n=1$.
 Without loss of generality, we assume that, for each $n\in \mathbb{N}$, $f^n$ satisfies (Hq') with the same Lipschitz constants.  \\

Obviously, ${f}^n$ has bounded first derivatives in $y$ and $z$ on the set $[0, T]\times \mathbb{R}^N\times [-{M}_y, {M}_y]\times [-{M}_z, {M}_z]$, for any ${M}_y$, ${M}_z>0$, because that $\hat{f}$ satisfies (Hq'). We now calculate the second derivative of ${f}^n$ in $y$ by\\
$$ \frac{\partial^2 {f}^{n}}{\partial y^2}=\int_{\mathbb{R}^2}
 \hat{f}\left({t}, x_1, x_2, \ldots, x_N, \tilde{y}, \tilde{z}\right)\frac{\partial^2 \rho_{n}}{\partial y^2}\left(y-\tilde{y}, z-\tilde{z} \right) d\tilde{y}d\tilde{z}.$$
Since on the set $[0, T]\times \mathbb{R}^N\times [-M_y-\frac{1}{n}, M_y+\frac{1}{n}]\times [-M_z-\frac{1}{n}, M_z+\frac{1}{n}]$,
$$
\left|\hat{f}\left( {t}, x_1, x_2, \ldots, x_N, {y}, {z}\right)\right|\leq M_0+L_y\bigg(M_y+1\bigg)+L_z\left(M_z+2\right)^2:=M_{yz},
$$
we have
$$
\bigg|\frac{\partial^2 {f}^{n}}{\partial y^2}\bigg|\leq C(n)M_{yz}\left|\left|\frac{\partial^2 \rho_{n}}{\partial y^2}\right|\right|_\infty,
$$
which means that $f^n$ has bounded second derivative in $y$ on the set $[0, T]\times \mathbb{R}^N\times [-{M}_y, {M}_y]\times [-{M}_z, {M}_z]$, for any ${M}_y$, ${M}_z>0$. Similar result can be obtained for the second derivative of ${f}^n$ in $z$.\\

From the results in the last section, we know that, for any $p\geq 2$, the GBSDE (\ref{po}) with the coefficient ${f}^n$ admits a solution $({Y}^n, {Z}^n, {K}^n)\in \mathfrak{G}^p_G(0, T)$. We now verify that $\{Y^n\}_{n\in \mathbb{N}}$ is a Cauchy sequence in $\mathcal{S}^p_G(0, T)$. For simplicity, we denote $f^n_t(y, z):=f^n(t, B_{t_1\wedge t}, B_{t_2\wedge t}-B_{t_1\wedge t}, \ldots, B_{t_{N}\wedge t}-B_{t_{{N-1}}\wedge t}, y, z)$ and $\hat{f}_t(y, z):=\hat{f}(t, B_{t_1\wedge t}, B_{t_2\wedge t}-B_{t_1\wedge t}, \ldots, B_{t_{N}\wedge t}-B_{t_{{N-1}}\wedge t}, y, z)$. Indeed, for $n$, $m\in \mathbb{N}$, $n\geq m$, $0\leq t\leq T$,
\begin{align}\label{recall1}
 \left|{f}^{n}_t\left({Y}^n_t,{Z}^n_t\right)- {f}^{m}_t\left({Y}^n_t,{Z}^n_t\right)\right|
&\leq |{f}^{n}_t\left({Y}^n_t,{Z}^n_t \right)- \hat{f}_t\left({Y}^n_t,{Z}^n_t \right)|+ |{f}^{m}_t\left({Y}^n_t,{Z}^n_t \right)- \hat{f}_t\left({Y}^n_t,{Z}^n_t\right)|\notag\\
 &\leq\int_{\mathbb{R}^2} \left|\hat{f}_t\left({Y}^n_t,{Z}^n_t \right)-\hat{f}_t\left({Y}^n_t-\tilde{y},{Z}^n_t-\tilde{z}\right)\right| \rho_{n}\left(\tilde{y}, \tilde{z}\right) d\tilde{y}d\tilde{z}\\
 &+\int_{\mathbb{R}^2} \left|\hat{f}_t\left({Y}^m_t,{Z}^m_t \right)-\hat{f}_t\left({Y}^m_t-\tilde{y},{Z}^m_t-\tilde{z}\right)\right| \rho_{m}\left(\tilde{y}, \tilde{z}\right) d\tilde{y}d\tilde{z}\notag\\
&\leq \frac{2}{m}(L_y+2L_z(M_z+1)),\notag
 \end{align}
 where $M_z$ is the bounded of all the solutions ${Z}^n$ defined in (\ref{solbou}), which is independent of $n$. Then, from Proposition \ref{esty}, we have $\{Y^n\}_{n\in \mathbb{N}}$ is a Cauchy sequence in $\mathcal{S}^p_G(0, T)$, so that there exists a $\hat{Y}\in \mathcal{S}^p_G(0, T)$, such that
 \begin{equation}\label{cony1}
 \mathbb{E}\left[\sup_{0\leq t\leq T}|Y^n_t-\hat{Y}_t|^p\right]\longrightarrow 0.
 \end{equation}

 Furthermore, from Proposition \ref{estz}, we have, for $1\leq p'/2<p$, $\{Z^n\}_{n\in \mathbb{N}}$  is a Cauchy sequence in $\mathcal{H}^{p'}_G(0, T)$. Then, there exists a $\hat{Z}\in \mathcal{H}^{p'}_G(0, T)$, such that
  \begin{equation}\label{conz1}
 \mathbb{E}\left[\left(\int^T_0|Z^n_t-\hat Z_t|^2dt\right)^{p'/2}\right]\longrightarrow 0.
 \end{equation}
 We define
 $$
 \hat{K}_t= \hat{Y}_t-\hat{Y}_0+\int^t_0 \hat{f}(s, B_{t_1\wedge s}, B_{t_2\wedge s}-B_{t_1\wedge s}, \ldots, B_{t_{N}\wedge s}-B_{t_{{N-1}}\wedge s}, \hat{Y}_s,\hat{Z}_s)d\langle B\rangle_s-\int^t_0 \hat Z_sdB_s.
 $$
 We proceed to prove that $\hat K$ is a decreasing $G$-martingale starting from $K_0=0$ and $\hat K_T\in L^{p'/2}_G(\Omega_T)$.
 it suffices to prove that
 $$
 \mathbb{E}\left[\left(\int^T_0\left|{f}^{n}_t\left({Y}^n_t,{Z}^n_t \right)-\hat{f}_t\left(\hat{Y}_t,\hat{Z}_t \right)\right|dt\right)^{p'/2}\right]\longrightarrow 0,
 $$
 which can be easily deduced by recalling (\ref{recall1}), (\ref{cony1}), (\ref{conz1}) and (Hq'). Thus, we obtain for $0\leq t\leq T$,
 $$
 \mathbb{E}\left[\left|K^n_t-\hat{K}_t\right|^{p'/2}\right]\longrightarrow 0,
 $$
from which have that, for $0\leq s\leq t\leq T$, $\mathbb{E}_s[\hat K_t]=\hat K_s$ by applying Proposition 2.5 in Hu et al. \cite{HJPS14a} and $\hat{K}_0=0$. Note that $p$ could be arbitrary large, which is ensured by the uniform boundedness of $Y^n$ and the $G$-BMO norm of $Z^n$. Therefore, for any $p\geq 2$, we could find a triple $(\hat{Y}, \hat{Z}, \hat{K})\in \mathfrak{G}^p_G(0, T)$, which is a solution to the GBSDE (\ref{po}) under Assumption \ref{apde} and \ref{cpde}.\\

{\bf Step 2}: Fix $N\in \mathbb{N}$ and the partition $\pi^N$ on $[0, T]$. We consider the GSDE (\ref{po}) with a generator $\bar{f}$ satisfying Assumption \ref{apde}.
In this step, we shall construct the solution to such a GBSDE by regularizing $\bar{f}$ in $t$ and $x_1, x_2\ldots, x_N$: for each $(y, z)\in \mathbb{R}^2$, $(t, x_1, x_2,\ldots, x_N)\in [0, T]\times\mathbb{R}^{N}$, we define
\begin{align*}
\hat{f}^{n}\left(t, x_1, x_2,\ \ldots, x_N,y,z \right)
:=\int_{\mathbb{R}^{N+1}} \bar{f}\left( t-\tilde{t}, x_1-\tilde{x}_1, x_2-\tilde{x}_2, \ldots, x_N-\tilde{x}_N,y,z\right)\\\times\rho_{n}\left( \tilde{t}, \tilde{x}_1, \tilde{x}_2, \ldots,\tilde{x}_N\right) d\tilde{t}d\tilde{x}_1d\tilde{x}_2\ldots d\tilde{x}_N,
\end{align*}
where $\rho_n$ is a positive smooth function such that its support is contained in a $\frac{1}{n}$-ball in $\mathbb{R}^{N+1}$ and $\int_{\mathbb{R}^{N+1}} \rho_n=1$. In addition, we define the extension of the function $\bar{f}$ on $\mathbb{R}_-$, i.e., if $t<0$, $\bar{f}(t, \cdot, \cdot, \ldots, \cdot, \cdot, \cdot):= \bar{f}(0, \cdot, \cdot, \ldots, \cdot, \cdot, \cdot)$. \\


Proceeding the same argument as in the last step, we can show that the first derivative of $\hat{f}^n$ in $t$, the first and second derivatives of $\hat{f}^n$ in $x_i$, $i=1, 2, \ldots, N$, are bounded on the set $[0, T]\times \mathbb{R}^N\times [-{M}_y, {M}_y]\times [-{M}_z, {M}_z]$, for any $M_y$, $M_z>0$. Therefore, recalling the result in the last step, we obtain that, for any $p\geq 2$, the GBSDE (\ref{po}) with the coefficient $\hat{f}^n$ admits a solution $(\hat{Y}^n, \hat{Z}^n, \hat{K}^n) \in \mathfrak{G}^p_G(0, T)$.\\

For $n$, $m\in \mathbb{N}$, $n\geq m$, $0\leq t\leq T$, by the definition of $\hat{f}^n$ and (Hc'),
\begin{align*}
|\hat{f}^n(t, B_{t_1\wedge t}, B_{t_2\wedge t}&-B_{t_1\wedge t}, \ldots, B_{t_{N}\wedge t}-B_{t_{{N-1}}\wedge t}, \hat{Y}^n_t, \hat{Z}^n_t)\\
&-
\hat{f}^m(t, B_{t_1\wedge t}, B_{t_2\wedge t}-B_{t_1\wedge t}, \ldots, B_{t_{N}\wedge t}-B_{t_{{N-1}}\wedge t}, \hat{Y}^n_t, \hat{Z}^n_t)|\\
&\leq (N+1)w^{\bar{f}}\left(\frac{1}{n}\right)\leq (N+1)w^{\bar{f}}\left(\frac{1}{m\wedge n}\right),
\end{align*}
from which we could deduce that $\{\hat{Y}^n\}_{n\in \mathbb{N}}$ is a Cauchy sequence, and similarly to (\ref{cony1}), there exists $\bar{Y}\in \mathcal{S}^p_G(0, T)$, such that
 \begin{equation*}\label{cony2}
 \mathbb{E}\left[\sup_{0\leq t\leq T}|\hat{Y}^n_t-\bar{Y}_t|^p\right]\longrightarrow 0.
 \end{equation*}
 We could conclude in a similar way as in the previous step that for any $p\geq 2$, there exists a triple $(\bar{Y}, \bar{Z}, \bar{K})\in \mathfrak{G}^p_G(0, T)$,
which solves the GBSDE (\ref{po}) under Assumption \ref{apde}.\\
%
%
%
%
%
%

To study the more general GBSDE (\ref{qgbsde}), we need the following assumption on the terminal value $\xi$.
\begin{assumption}\label{assxi}
$\xi\in L^\infty_G$. 
\end{assumption}
We are now ready to introduce the main result of this paper.
\begin{theorem}
Consider the GBSDE (\ref{qgbsde}) satisfying Assumption \ref{assbig} and \ref{assxi}. It admits at least a solution $({Y}, {Z}, {K})\in \mathfrak{G}^2_G(0, T)$.
\end{theorem}
\noindent {\bf Proof:}
%
Without loss of generality, we assume $\xi$ is approximated by the following sequence
$$
\xi^n:=\varphi^n(B_{t^n_1}, B_{t^n_2}-B_{t^n_1}, \ldots, B_{t^n_{N(n)}}-B_{t^n_{N(n)-1}})\in Lip(\Omega_T),
$$
where for each $n\in \mathbb{N}$, $0=t^n_0\leq t^n_1\leq \ldots\leq t^n_{N(n)}=T$, $\mu^n:=\max_{k=1, 2, \ldots, N(n)}|t^n_k-t^n_{k-1}|\leq 1/2^{n}$. Assume moreover that for each $n\geq m$, $\{0=t^m_0, t^m_1, t^m_2, \ldots, t^m_{N(m)}=T\}=:\pi^m \subset \pi^n:=\{0=t^n_0, t^n_1, t^n_2, \ldots, t^n_{N(n)}=T\}$.\\

Fix $n\in \mathbb{N}$. We construct the function $\bar f^n$ in terms of $h$ by discretization. For simplicity of notation, we omit the superscript $n$ for $t^n_k$, $k=0, 1, \ldots, {N(n)}$.
%
Denote by ${\bf x}(n)$ the vector $(x_1, x_2, \ldots, x_{N(n)})\in \mathbb{R}^{N(n)}$. Let $t\in [t_{k-1}, t_{k}]$, where $k=1, 2, \ldots \leq N(n)$. We define by the following procedures a piecewisely linear path stopped at time $t$ in terms of ${\bf x}(n)$, noted by  $\omega^{{\bf x}(n), t}$.
\begin{itemize}
\item $\omega^{{\bf x}(n), t}_{t_0}=0$;
\item $\omega^{{\bf x}(n), t}_{t_1}=x_1$;
\item $\omega^{{\bf x}(n), t}_{t_2}=x_1+x_2$;
\item $\ldots$;
\item $\omega^{{\bf x}(n), t}_{t_{k-1}}=\sum^{k-1}_{i=1} x_i$;
\item $\omega^{{\bf x}(n), t}_t=\omega^{{\bf x}(n), t}_{t_{k}}=\omega^{{\bf x}(n), t}_{t_{k+1}}=\omega^{{\bf x}(n), t}_{t_{k+2}}=\ldots=\omega^{{\bf x}(n), t}_{t_{N(n)}}=\sum^{k}_{i=1} x_i$;
\item $\omega^{{\bf x}(n), t}$ is a linear function in $t$ on $[t_{i-1}, t_{i}]$, for $i=1, 2, \ldots, k$. It is also linear on $[t_{k}, t]$ and takes a constant value on $[t, T]$.
\end{itemize}


\noindent Define $$
\bar{f}^n(t, x_1, x_2, \ldots, x_{N(n)}, y, z):= h(t, \omega^{{\bf x}(n), t}, y, z).
$$
We can verify that for each $t\in [0, T]$, ${\bf x}^1(n):=(x^1_1, x^1_2, \ldots, x^1_{N(n)})$, ${\bf x}^2(n):=(x^2_1, x^2_2, \ldots, x^2_{N(n)})$ and $(y, z)\in \mathbb{R}^2$,
\begin{align*}
\bigg|\bar{f}^n\left(t, x^1_1, x^1_2, \ldots, x^1_{N(n)}, y, z\right)&-\bar{f}^n\left(t, x^2_1, x^2_2, \ldots, x^2_{N(n)}, y, z\right)\bigg|\\
&\leq w^h\left(\left|\left|\omega^{{\bf x}^1(n), t}-\omega^{{\bf x}^2(n), t}\right|\right|_\infty\right)\leq w^h\left( \sum_{k=1, 2, \ldots, N(n)}|x^1_k-x^2_k|\right),
\end{align*}
which implies that $\bar{f}$ is uniformly continuous with modulus $w^h$ in $x_1$, $x_2, \ldots, x_{N{(n)}}$, where the modulus is independent of $y$ and $z$. From the result in Step 2, we know that, for any $p\geq 2$, the GBSDE (\ref{po}) with the parameters $(\bar{f}^n, \xi^n)$ admits a solution $(\bar{Y}^n, \bar{Z}^n, \bar{K}^n)\in \mathfrak{G}^p_G(0, T)$.\\

%
%
%
For $n$, $m\in \mathbb{N}$, $n\geq m$, we denote
\begin{align*}
\eta^{n, m}&:=\int^T_0\bar{f}^n(t, B_{t_1\wedge t}, B_{t_2\wedge t}-B_{t_1\wedge t}, \ldots, B_{t_{N(n)}\wedge t}-B_{t_{{N(n)-1}}\wedge t}, \bar{Y}^n_t, \bar{Z}^n_t)d\langle B\rangle_t\\
&-\int^T_0\bar{f}^m(t, B_{t_1\wedge t}, B_{t_2\wedge t}-B_{t_1\wedge t}, \ldots, B_{t_{N(m)}\wedge t}-B_{t_{{N(m)-1}}\wedge t}, \bar{Y}^n_t,\bar{Z}^n_t)d\langle B\rangle_t.
\end{align*}
By Proposition \ref{esty}, we have
\begin{align*}
| Y^n_{t}-Y^m_{t}| \leq C\left( \Vert \xi ^{n}-\xi ^{m}\Vert_{L^\infty_G}+\tilde{\mathbb{E}}_t\left[ \eta^{n, m}\right] \right)
\leq C\left( \Vert \xi ^{1}-\xi ^{2}\Vert_{L^\infty_G}+{\mathbb{E}}_t\left[ \frac{\mathcal{E}(b^{n, m})_T}{\mathcal{E}(b^{n, m})_t}\eta^{n, m}\right] \right),
\end{align*}
where $b^{n,m}$ is defined by (\ref{poe}) in terms of $\bar{Z}^n$ and $\bar{Z}^m$. Fortunately, its $G$-BMO norm is dominated by a constant $C_b$ independent of $n$ and $m$. Thus, we could find a uniform order $q>1$ and a uniform constant $C_q$ for the reverse H\"older inequality for $\mathcal{E}(b^{n, m})$, for all $n, m\in \mathbb{N}$. Since the function $\Phi$ is decreasing (see Theorem 3.1 in \cite{Kaz94}), we assume without loss of generality that $q<2$. Then, we obtain
\begin{align*}
| Y^n_{t}-Y^m_{t}|
\leq C\left( \Vert \xi ^{1}-\xi ^{2}\Vert_{L^\infty_G} + C^{1/q}_q {\mathbb{E}}_t \left[(\eta^{n, m})^{p}\right]^{1/p} \right),
\end{align*}
where $1/q+1/{p}=1$. Furthermore, applying Theorem 2.8 in \cite{HJPS14a}, we have, for $2<p<p''$ and $1<\gamma<p''/p$,
\begin{align*}
\mathbb{E}\left[\sup_{0\leq t\leq T}| \bar{Y}^n_{t}-\bar{Y}^m_{t}|^p\right]
&\leq C\left(\Vert \xi ^{1}-\xi ^{2}\Vert^p_{L^\infty_G} + \mathbb{E}\left[\sup_{0\leq t\leq T}{\mathbb{E}}_t \left[(\eta^{n, m})^{p}\right]\right] \right)\\
&\leq C\left(\Vert \xi ^{1}-\xi ^{2}\Vert^p_{L^\infty_G} + \mathbb{E}\left[(\eta^{n, m})^{p''}\right]^{p/p''}+ \mathbb{E}\left[(\eta^{n, m})^{p''}\right]^{1/\gamma}\right),
\end{align*}
where the constant $C$ varies from line to line, nevertheless, is independent of $n$ and $m$.\\
%

To verify that $\{\bar{Y}^n\}_{n\in \mathbb{N}}$ is a Cauchy sequence in $\mathcal{S}^p_G(0, T)$, it suffices to have
\begin{equation}\label{coneta}\mathbb{E}\left[(\eta^{n, m})^{p''}\right]\longrightarrow 0,\quad \mbox{as}\quad m, n\rightarrow \infty.\end{equation}
For simplicity, note $B^t:=B_{\cdot\wedge t}$. For each $\omega\in \Omega$, we define $${\bf x}^{B^t(\omega)}(n):= (B_{t_1}(\omega), B_{t_2}(\omega)-B_{t_1}(\omega), \ldots, B_{t}(\omega)-B_{t_{k-1}}(\omega), \underbrace{0, \ldots, 0}_{N(n)-k}),\ t\in [t_{k-1}, t_{k}]$$
and a mapping ${B}^{n, t}: \mathcal{C}(0, T)\rightarrow \mathcal{C}(0, T)$,
$${B}^{n, t}(\omega):= \omega^{{\bf x}^{B^t(\omega)}(n), t},$$
Then,
\begin{align*}
\eta^{n, m}&=\int^T_0\left|\bar{f}^n(t, {\bf x}^{B^t(\cdot)}(n), \bar{Y}^n_t, \bar{Z}^n_t)-\bar{f}^m(t, {\bf x}^{B^t(\cdot)}(m), \bar{Y}^n_t, \bar{Z}^n_t)\right|d\langle B\rangle_t\\
&=\int^T_0\left|h(t, {B}^{n, t}(\cdot), \bar{Y}^n_t, \bar{Z}^n_t)-h(t, {B}^{m, t}(\cdot), \bar{Y}^n_t, \bar{Z}^n_t)\right|d\langle B\rangle_t\\
&\leq \int^T_0 w^h\left(\left|\left|B^{n, t}(\cdot)-B^{m, t}(\cdot)\right|\right|_{\infty}\right)dt\\
&\leq \int^T_0 w^h\left(\left|\left|B^{n, t}(\cdot)-B^t(\cdot)\right|\right|_{\infty}\right)dt+  \int^T_0 w^h\left(\left|\left|B^{m, t}(\cdot)-B^t(\cdot)\right|\right|_{\infty}\right)dt,
\end{align*}
where the last inequality is from the sub-additivity of $w^h$.\\

Now, our aim is to prove that
$$
\mathbb{E}\left[\int^T_0 w^h\left(\left|\left|B^{n, t}(\cdot)-B^t(\cdot)\right|\right|_{\infty}\right)dt\right]\longrightarrow 0,\quad\mbox{as}\quad n\rightarrow \infty.
$$
Then, due to the boundedness of $w^h$, we have
\begin{align}\label{conmid3}\mathbb{E}\left[(\eta^{n, m})^{p''}\right]&\leq C_{p''}\bigg(\mathbb{E}\left[\left(\int^T_0 w^h\left(\left|\left|B^{n, t}(\cdot)-B^t(\cdot)\right|\right|_{\infty}\right)dt\right)^{p''}\right]\\
&\quad\quad\quad+\mathbb{E}\left[\left(\int^T_0 w^h\left(\left|\left|B^{m, t}(\cdot)-B^t(\cdot)\right|\right|_{\infty}\right)dt\right)^{p''}\right]\bigg)\longrightarrow 0, \notag
\quad\mbox{as}\quad m, n\rightarrow \infty.\end{align}
Indeed, since $w^h$ is a concave function, by Lemma 2.12 in Bai and Lin \cite{BL14}, we have
\begin{align}
\mathbb{E}\left[\int^T_0 w^h\left(\left|\left|B^{n, t}(\cdot)-B^t(\cdot)\right|\right|_{\infty}\right)dt\right]&\leq
w^h\left(\mathbb{E}\left[\int^T_0 \left|\left|B^{n, t}(\cdot)-B^t(\cdot)\right|\right|_{\infty}dt\right]\right)\label{conmid}\\
&\leq w^h\left(\int^T_0 \mathbb{E}\left[\left|\left|B^{n, t}(\cdot)-B^t(\cdot)\right|\right|_{\infty}\right]dt\right).\notag
\end{align}
Then, to eventually prove (\ref{coneta}), it suffices to have
\begin{equation}
\mathbb{E}\left[\left|\left|B^{n, t}(\cdot)-B^t(\cdot)\right|\right|_{\infty}\right]
\longrightarrow 0,\quad\mbox{as}\quad n\rightarrow \infty, \label{conmid2}\end{equation}
and then we could apply the Lebesgue dominated convergence theorem to deduce the convergence of (\ref{conmid}), which yields (\ref{conmid3}).
For $t\in [t_{k-1}, t_{k}]$, $k=1, 2, \ldots, N(n)$, we calculate (\ref{conmid2}) and obtain
\begin{align*}
&\ \ \ \ \Vert B^{n, t}(\omega)-B^t(\omega)\Vert_\infty\\
&\leq \max_{i=1, 2, \ldots, k-1} \left(\sup_{s\in [t_{i-1}, t_{i})}B_s(\omega)-\inf_{s\in [t_{i-1}, t_{i})}B_s(\omega) \right)\vee \left(\sup_{s\in [t_{k-1}, t]}B_s(\omega)-\inf_{s\in [t_{k-1}, t]}B_s(\omega) \right)\\
&\leq 2\left(\max_{i=1, 2, \ldots, k-1} \left(\sup_{s\in [t_{i-1}, t_{i})}\left|B_s(\omega)-B_{t_{i-1}}(\omega)\right|\right)\vee\sup_{s\in [t_{k-1}, t]}\left|B_s(\omega) -B_{t_{k-1}}(\omega)\right|\right).
\end{align*}
Thus, for $\alpha>2$,
$$
\mathbb{E}[\Vert B^{n, s}(\cdot)-B^s(\cdot)\Vert^\alpha_\infty]\leq C_{\overline{\sigma}} 2^{-n(\frac{\alpha}{2}-1)}\longrightarrow 0, \quad\mbox{as}\quad{n\rightarrow\infty},
$$
which implies (\ref{conmid2}).\\

We conclude from (\ref{coneta}) that, for $p>2$, there exists a process $Y\in \mathcal{S}^p_G(0, T)$ such that
\begin{align*}
\mathbb{E}\left[\sup_{0\leq t\leq T}| \bar{Y}^n_{t}-{Y}_{t}|^p\right]\longrightarrow 0, \quad\mbox{as}\quad n\rightarrow \infty.\end{align*}
Then, proceeding the same argument as in Step 1, we obtain, for $2\leq p'/2<p$, there exists $Z\in\mathcal{H}^{p'}_G(0, T)$, such that
\begin{equation*}
 \mathbb{E}\left[\left(\int^T_0|\bar{Z}^n_t-Z_t|^2dt\right)^{p'/2}\right]\longrightarrow 0, \quad\mbox{as}\quad n\rightarrow \infty,
 \end{equation*}
 and
 $$
 \mathbb{E}\left[\left(\int^T_0\left|\bar{f}^{n}\left(t, {\bf x}^{B^t(\cdot)}(n), \bar{Y}^n_t,\bar{Z}^n_t \right)-{h}\left(t, B_{\cdot\wedge t}(\cdot), {Y}_t,{Z}_t \right)\right|dt\right)^{p'/2}\right]\longrightarrow 0,\quad\mbox{as}\quad n\rightarrow \infty.
 $$
 Thus, there exists process $K$ with $K_0=0$ and $K_T\in L^{p'/2}_G(\Omega_T)$, and due to the fact that 
 for $0\leq t\leq T$,
 $$
 \mathbb{E}\left[\left|\bar{K}^n_t-{K}_t\right|^{p'/2}\right]\longrightarrow 0,  \quad\mbox{as}\quad n\rightarrow \infty,
 $$
 which implies that $K$ is a decreasing $G$-martingale. 
 In conclusion, since $p$ and $p'$ can be arbitrarily large, we could finally find  the triple $({Y}, {Z}, {K})\in \mathfrak{G}^2_G(0, T)$ which solves (\ref{qgbsde}). We complete the proof. \hfill$\square$\\
 
 \noindent {\bf Acknowledgement:} Y. Hu and A. Soumana Hima are grateful for partial financial support from the Lebesgue Center of Mathematics (``Investissements d'avenir" Program) under grant ANR-11-LABX-0020-01 and  from the Agence Nationale de la Recherche (ANR) under grant ANR-15-CE05-0024-02. Yiqing Lin  gratefully acknowledges financial support from the European Research Council (ERC) under grant FA506041 and under grant 321111. We thank Falei Wang for helpful suggestions. 


\bibliography{Ref}
\bibliographystyle{plain}

\end{document}